\RequirePackage{fix-cm}
\documentclass[smallextended]{svjour3}       % onecolumn (second format)
\smartqed  % flush right qed marks, e.g. at end of proof
\usepackage{graphicx}

\usepackage[left=2.5cm,right=2.5cm,top=1cm,bottom=1cm,includeheadfoot]{geometry}

\usepackage{csquotes}
\usepackage[english]{babel}
\usepackage{amsmath, amssymb, mathtools, amsfonts,nicefrac}
\usepackage[mathscr]{eucal}
\usepackage{booktabs}

\usepackage{mathabx}
\usepackage{hyperref}
\hypersetup{hidelinks,colorlinks,linkcolor=blue,urlcolor=black}

\smartqed  

\newenvironment{proof}{\par\noindent\textbf{Proof.}}{\hfill$\square$\par}

\makeatletter
\def\cl@chapter{\@elt {theorem}}
\makeatother
\makeatletter
\let\cl@chapter\undefined
\makeatother
\usepackage[nameinlink]{cleveref}

\crefname{asmpt}{\textup{Assumption}}{\textup{Assumptions}}
\creflabelformat{asmpt}{#2\textup{#1}#3} 
\crefname{theorem}{\textup{Theorem}}{\textup{Theorems}}
\creflabelformat{theorem}{#2\textup{#1}#3} 
\crefname{lemma}{\textup{Lemma}}{\textup{Lemmata}}
\creflabelformat{lemma}{#2\textup{#1}#3} 
\crefname{proposition}{\textup{Proposition}}{\textup{Propositions}}
\creflabelformat{proposition}{#2\textup{#1}#3} 
\crefname{cor}{\textup{Corollary}}{\textup{Corollaries}}
\creflabelformat{cor}{#2\textup{#1}#3} 
\crefname{definition}{\textup{Definition}}{\textup{Definitions}}
\creflabelformat{definition}{#2\textup{#1}#3} 
\crefname{remark}{\textup{Remark}}{\textup{Remarks}}
\creflabelformat{remark}{#2\textup{#1}#3} 
\crefname{figure}{\textup{Figure}}{\textup{Figures}}
\creflabelformat{figure}{#2\textup{#1}#3} 
\crefname{table}{\textup{Table}}{\textup{Tables}}
\creflabelformat{table}{#2\textup{#1}#3} 
\crefname{algorithm}{\textup{Algorithm}}{\textup{Algorithms}}
\creflabelformat{algorithm}{#2\textup{#1}#3} 
\crefname{section}{\textup{Section}}{\textup{Sections}}
\creflabelformat{section}{#2\textup{#1}#3}

\usepackage{subcaption}
\usepackage[draft,todonotes={textsize=scriptsize}]{changes}
\makeatletter
\let\@citexOld\@citex
\def\@citex[#1]#2{\textup{\@citexOld[#1]{#2}}}
\makeatother

\usepackage{color}

\usepackage{algorithm} %provides pseudocode framework
\usepackage{algpseudocode} %\usepackage{algorithmic} 

\usepackage{enumitem}

\numberwithin{algorithm}{section}
\numberwithin{equation}{section}
\numberwithin{figure}{section}
\numberwithin{table}{section}
\numberwithin{table}{section}
\numberwithin{remark}{section}
\numberwithin{equation}{section}
\usepackage[toc]{appendix}

\counterwithin{rem}{section}

\newtheorem{asmpt}{Assumption}

\newcommand{\sC}{{\mathscr C}}
\newcommand{\sD}{{\mathscr D}}

\newcommand{\sJ}{{\mathscr{J}}}

\newcommand{\sP}{{\mathscr{P}}}
\newcommand{{\sQ}}{{\mathscr{Q}}}

\newcommand{\sT}{{\mathscr T}}

\newcommand{\sPad}{{\sP_{\mathrm{ad}}}}

\newcommand{\Nb}{\mathbb{N}}
\newcommand{\Rb}{\mathbb{R}}

\newcommand{\bb}{{\boldsymbol{b}}}

\newcommand{\bp}{{\boldsymbol{p}}}
\newcommand{\bq}{{\boldsymbol{q}}}

\newcommand{\bu}{{\boldsymbol{u}}}
 
\newcommand{\bz}{{\boldsymbol{z}}}

\newcommand{\tildeq}{{\widetilde{q}}}
\newcommand{\tildebq}{{\widetilde{\boldsymbol{q}}}}

\newcommand{\rmopt}{{\mathrm{opt}}}
\newcommand{\rmref}{{\mathrm{ref}}}

\newcommand{\Ttrans}{{T_{\mathrm{trans}}}}
\newcommand{\Ttot}{T_{\mathrm{tot}}}

\def\scalingThreeFigures{0.27}

\newcommand{\revisionComment}[1]{\textcolor{black}{#1}} %teal
\newcommand{\SecondRevisionComment}[1]{\textcolor{black}{#1}} %red

\allowdisplaybreaks

\begin{document}
	
	\title{Reconstructing the system coefficients for coupled harmonic oscillators}
	
	% Authors: full names plus addresses.
	\author{
        Jan Bartsch \and 
        Ahmed A. Barakat \and 
        Simon Buchwald \and 
        Gabriele Ciaramella \and ´
        Stefan Volkwein \and 
        Eva M. Weig
	}

 \institute{
        Jan Bartsch (corresponding author) \at 
        Department of Mathematics, University of Konstanz \\
        Universitätsstraße 10, 78464 Konstanz, Germany \\
        \email{jan.bartsch@uni-konstanz.de} \\
        ORCID: \url{https://orcid.org/0000-0002-8011-7422}
   \and 
        Ahmed A. Barakat \at 
        TUM School of Computation, Information and Technology, Technical University of Munich \\
        Hans-Piloty-Str. 1(5901)/III, 85748 Garching b. München, Germany 
        \at
        Design and Production Engineering, Faculty of Engineering, Ain Shams University\\
        El-Sarayat St. 1, Abbaseya, 11517 Cairo, Egypt.\\
        \email{ahmed.barakat@tum.de} \\
        ORCID: \url{https://orcid.org/0000-0003-2197-1124}
	\and 
        Simon Buchwald \at 
        Department of Mathematics, University of Konstanz \\
        Universitätsstraße 10, 78464 Konstanz, Germany \\
        \email{simon.buchwald@uni-konstanz.de} \\
        ORCID: \url{https://orcid.org/0009-0004-2350-4399}
    \and 
        Gabriele Ciaramella \at 
        MOX Lab, Dipartimento di Matematica, Politecnico di Milano\\
        Piazza Leonardo da Vinci 32, 20133, Milano, Italy \\
        Member of GNCS Indam group \\
        \email{gabriele.ciaramella@polimi.it} \\
        ORCID: \url{https://orcid.org/0000-0002-5877-4426}
    \and 
        Stefan Volkwein \at 
        Department of Mathematics, University of Konstanz \\
        Universitätsstraße 10, 78464 Konstanz, Germany \\
        \email{stefan.volkwein@uni-konstanz.de} \\
        ORCID: \url{https://orcid.org/0000-0002-1930-1773}
    \and 
        Eva Weig \at 
        TUM School of Computation, Information and Technology, Technical University of Munich \\
        Hans-Piloty-Str. 1, 85748 Garching b. München, Germany \\
        TUM Zentrum für QuantumEngineering (ZQE), Am Coulombwall 3a, 85748 Garching b. München, Germany\\
        Munich Center for Quantum Science and Technology (MCQST), Schellingstr. 4, 80799 München, Germany
        \email{eva.weig@tum.de} \\
        ORCID: \url{https://orcid.org/0000-0003-4294-8601}
	}
	
	\date{Received: date / Accepted: date}
	\maketitle

\begin{abstract}
    Physical models often contain unknown functions and relations.
    In order to gain more insights into the nature of physical processes, these unknown functions have to be identified or reconstructed.
    Mathematically, we can formulate this research question within the framework of inverse problems. 
    In this work, we consider optimization techniques to solve the inverse problem using Tikhonov regularization and data from laboratory experiments. We propose an iterative strategy that eliminates the need for \SecondRevisionComment{further} laboratory experiments. 
    Our method is applied to identify the coupling and damping coefficients in a system of oscillators, ensuring an efficient and experiment-free approach. 
    We present our results and compare them with those obtained from an alternative, purely experimental approach. 
    By employing our proposed strategy, we demonstrate a significant reduction in the number of laboratory experiments required.
\end{abstract}



%%%%%%%%%%%%%%%%%%%%%%%%%%%%%%%%%%%%%%
\section{Introduction}
\label{sec:intro}
%%%%%%%%%%%%%%%%%%%%%%%%%%%%%%%%%%%%%%

The modeling of physical phenomena is crucial for understanding natural processes, allowing to predict behaviors and outcomes under various conditions. 
These models often rely on a set of coefficients that describe interactions and dynamics within the physical model. 
The accurate identification of these coefficients is essential for constructing reliable models that can simulate real-world behavior.

However, physical models often contain unknown functions and relations, posing significant challenges. 
These unknowns can stem from complex interactions within the model that are not directly observable. 
In order to gain more insights into the nature of physical processes, these unknown functions must be identified. 
This identification procedure involves adjusting the model parameters to ensure that the predictions of the model align with experimental observations. 
Consistency between these predictions and laboratory experiments is critical for validating both the accuracy of the reconstructed coefficients and the overall model.

Mathematically, we can formulate this research question within the framework of inverse problems; see, e.g., \cite{Kirsch2021IntroInverseProblems}.
Inverse problems involve determining the causal factors from observed effects, and are often ill-posed, meaning that solutions may not exist, be unique, or depend continuously on the data. 
Regularization techniques, such as Tikhonov regularization, are employed to stabilize the solution of these problems by incorporating prior knowledge or additional constraints \cite{Kirsch2021IntroInverseProblems}.

In this work, we apply this framework to a laboratory experiment and work with \SecondRevisionComment{real-world measurement} data generated in this experiment.

We investigate a system of coupled linear driven and damped oscillators in which the damping and coupling coefficients are unknown. 
The linear oscillators correspond to the two fundamental vibrational modes of a nanomechanical resonator visualized in \cref{fig:VibrationModes}. 
It consists of a long and thin nanostring which is a few hundred nanometers wide and thick and several tens of micrometers long and which vibrates at megahertz frequencies.
The coupling between the modes is mediated by an electrical control field. 
Such nanoelectromechanical systems are employed, e.g., for sensing applications~\cite{BachtoldMoserDykman2022ReviewMesoNanomechanical,barzanjeh_optomechanics_2022,guo2023mode}, such that a detailed understanding and modeling of its  parameters is of crucial interest. 
Moreover, since the approach in this work is essentially mathematical, it can certainly be adjusted and extended to investigate the coupling interaction in other two-mode, or two-level, physical systems.
In addition, similar signatures can be found, e.g., in cavity optomechanical systems or in circuit quantum electrodynamics \cite{aspelmeyer_cavity_2014,clerk_hybrid_2020}.

While the uncoupled system is very well understood \cite{GajoRastelliWeig2019NonlinearCouplingResonator,OchsWeig2022FrequencyNonlinNanoMode,OchsWeig2021NonlinResponseBrokenSymmetrie}, this is not the case for the coupled one.
The identification of coefficients in the coupled oscillator and similar problems have been addressed in the literature \cite{KralemannCimponeriu2008PhaseDynamicsCoupledOsc,PanoggioCiocanelLazarus2019ReconOsc}.
Our approach leverages both simulation and experimental data to enhance the accuracy of the coefficient reconstruction. 
It serves as a proof-of-principle demonstration and allows to be applied to more complex situations, where a conventional identification of coefficients via calibration measurements is no longer feasible.

Since the data from our experiments inherently include the effects of coupling, we cannot measure the eigenfrequencies of the system and the coupling coefficient separately. 
On the contrary, our experiments only allow us to have access to a hybridized quantity that contains both the eigenfrequencies and coupling coefficient.
Therefore, these parameters and coefficients must be inferred through indirect methods that reconstruct them from measurement data while accounting for the coupling coefficient and eigenfrequency separately. 
\SecondRevisionComment{Furthermore, we need to consider a scaling factor between the simulation and the experiments that is also depending on the unknown parameters. This leads to a highly nonlinear optimization problem. 
To tackle this, we need to carefully state the optimization problem and formulate a new algorithm that is able to solve this problem.
Moreover, we take into account that the system under investigation undergoes changes during  the conducting of experiments. More specifically, the eigenfrequencies drift over the time of the experiments.
The detailed steps in our approach to model these issues  are not common in the literature but  fundamental to problems of this form when working with real experiments.}

\SecondRevisionComment{With regard to the challenges outlined above in reconstructing the unknown coupling coefficient, it becomes clear that applying the reconstruction method to the physical problem introduced additional constraints that are not accounted for in the theory and therefore have to be incorporated into the model. Identifying these constraints enabled the successful integration of experiments and optimization and hence represents the core contribution of the present work. This contribution is significant because it addresses a problem in physics application for which existing methods often are inadequate or demand greater experimental cost~\cite{BarakatWeig2024CouplingStrength}. In addition, we emphasize the importance of testing optimization methods on real physical problems, as this provides the practical framework that is otherwise missing for their effective application.}

Our proposed strategy involves an iterative optimization technique using Tikhonov regularization to solve the inverse problem \cite{HaberAscherOldenburg2000OptimizationNonlinInverse,KunischZou1998IterativeRegularizationInverseProblems}. 
\revisionComment{We remark again that the parameter identification problem  and the mathematical optimization method are not novel.
The novelty is to build a concrete optimization problem that corresponds exactly to the laboratory experiment motivated by fundamental physics.
Usually, the mathematical models are only loosely related the laboratory settings and miss important difficulties (e.g. hypridized model in experiments, scaling factor, drifting eigenfrequencies).
In this work, we address all issues arising in the experimental setup. 
In particular, we consider the hybridization of the two modes and the need for reconstructing the transformation matrix $\mathscr{T_\theta}$ (see \cref{sec:Setting}).
Furthermore, our work shows that with our concrete model, which has a structure that is well-known in the experimentalist community, it is sufficient to apply standard optimization tools for  the reconstruction of parameters.
}

To generate measurement data, we perform several experiments in the laboratory with different external excitation forces.

An additional important new contribution of our work is that the presented approach allows us to decrease the number of required laboratory experiments, reducing experimental costs and time. 
In particular, our approach needs fewer experiments than the one described in \cite{BarakatWeig2024CouplingStrength}.
Furthermore, our experiments can be conducted easily and fast.
Additionally, after fixing the physical model in which parameters should be identified, no additional physical concepts need to be taken into account.
We compare our results to those obtained from a purely experimental approach that needs a larger number of experiments that are also more difficult to conduct \cite{BarakatWeig2024CouplingStrength}.
By this comparison, we demonstrate the effectiveness of our method in reconstructing the system coefficients for coupled harmonic oscillators.

Moreover, our approach can be extended to different (also nonlinear) physical systems as long as they can be sufficiently well described by ordinary differential equations and the corresponding experiments can be conducted quickly.
Also in view of this extension, our current work serves as a proof-of-principle.

The article is organized as follows. 
In \cref{sec:Setting}, we introduce our specific framework and provide details on the transformation between the measurable model to which we have access in the experiments and the inaccessible one to which we want to identify the coefficients.
\SecondRevisionComment{Moreover, we discuss the specific external driving used in the experiments, how the conducting of experiments lead to slight changes in the eigenfrequencies, and how we deal with this in the problem formulation.}
We describe the physical background and laboratory setup in \cref{sec:PhysicalSystem} before we explain in \cref{sec:FRF_calibration} how we scale the simulation data in order to \SecondRevisionComment{quantitatively} compare them with the experimentally measured physical data. 
%\EWnote{Question EW: ones - or physical data}
\SecondRevisionComment{Since this scaling also depends on the unknown parameters, it will also be included in the algorithm solving the inverse problem.}
In \cref{sec:ParamIdentProblem}, we formulate the (regularized) inverse problem and present an algorithm to solve it.
\SecondRevisionComment{This is a vital part of the manuscript, since the correct problem formulation is crucial  given its physical nature. In particular, identifying precisely in which parts of the models unknown parameters are included is not trivial in this case.}
Afterward, we discuss in \cref{sec:Solving_ODEs} how we solve the system of equations numerically \SecondRevisionComment{with a procedure that can be extended to other ODEs.}
Finally, in \cref{sec:NumExp}, we perform the final identification of the system coefficients by presenting our results and shortly compare it to a purely experimental approach in \cref{sec:Comparison}.
A section of conclusions completes this work.

\medskip
\textbf{Notation.} 
For a vector $\bp \in \Rb^n$, we denote by $\|\bp\|_2$ the standard Euclidean norm: 
$\|\bp\|_2 = \sqrt{\sum_{i=1}^n p_i^2}$.
For a vector $\bp$ or a matrix $A$, we denote by '$\bp^{\top}$' or '$A^{\top}$' their transpose.
Whenever a vector or matrix is compared or divided by a vector or matrix, respectively, these operations have to be understood componentwise.

%%%%%%%%%%%%%%%%%%%%%%%%%%%%%%%%%%%%%%
\section{System of coupled oscillators}	
\label{sec:Setting}
%%%%%%%%%%%%%%%%%%%%%%%%%%%%%%%%%%%%%%

In this section, we introduce the setting of the considered parameter estimation problem that will be formulated in detail in \cref{sec:ParamIdentProblem}. 
We consider a vibrating nanostring that is modeled by a two-dimensional damped and driven harmonic oscillator.
First, in \cref{sec:StandardModel}, the model is formulated in its standard form in which the coupling appears explicitly. 
However, the observation during measurement considers the system in transformed coordinates, where the pure mechanical eigenmodes of the oscillator are hybridized, or in a mathematical sense, where the stiffness matrix is diagonalized. 
Because of the hybridization of the physical coordinates, we need to transform the physical coordinates into hybridized modal coordinates, where the coupling remains only implicitly effective.
We present the transformation of coordinates in detail in \cref{sec:Transformation}.
The structure of the driving forces that we use in the experiments is explained in \cref{sec:DrivingForces}.

%%%%%%%%%%%%%%%%%%%%%%%%%%%%%%%%%%%%
\subsection{Standard model}
\label{sec:StandardModel}
%%%%%%%%%%%%%%%%%%%%%%%%%%%%%%%%%%%

We consider a coupled system of linear second-order ordinary differential equations (ODEs) in the physical coordinates $q_1$ and $q_2$.
%that represent the eigenmodes of the system.
The evolution of $\bq\coloneqq (q_1,q_2)^\top$ over a finite-time horizon $\Ttot>0$ is then given by the coupled linear second-order system 
\begin{align}
   \binom{\ddot{q}_1(t)}{\ddot{q}_2(t)}
    +
    \underbrace{\begin{pmatrix}
            2\pi d_1  & 0
            \\
            0 & 2\pi d_2 
    \end{pmatrix}}_{\eqqcolon \mathscr{D}}
    \binom{\dot{{q}}_1(t)}{\dot{{q}}_2(t)}
    + \underbrace{\begin{pmatrix}
            (2\pi f_{1})^2 & -(2\pi \lambda)^2
            \\
            -(2\pi \lambda)^2 & (2\pi f_{2})^2
    \end{pmatrix}}_{\eqqcolon \sC}
    \binom{q_1(t)}{q_2(t)}
    = \underbrace{\binom{b_1(t)}{b_2(t)}}_{\eqqcolon \bb},
    \qquad t \in [0,\Ttot]
    ,
     \label{eq:original_system}
\end{align}	
with initial condition $q_1(0) = 0 = q_2(0)$ and $\dot{q}_1(0) = 0 = \dot{q}_2(0)$.
Since the driving $\bb$ is assumed to be smooth, the solution $\bq$ to \eqref{eq:original_system} is unique and smooth.
For further analysis of such systems, see, e.g., \cite{Borzi2020DGLBook,Singiresu1995mechanicalVibration}.

As the observation of the system, we consider its Fourier transformation.
To avoid transient effects within the Fourier transform, we will discard the transient phase $[0,\Ttrans]$, $0\ll\Ttrans\ll\Ttot$, and only consider the system in equilibrium for a fixed amount of time $T$, i.e., in the interval $[\Ttrans,\Ttrans+T]$, where we have defined $T\coloneqq \Ttot-\Ttrans$.
Hence, the initial condition for \eqref{eq:original_system} is just given for the sake of completeness and the fact that we need to specify it for the numerical simulations.

Notice that the constants $d_1$ and $d_2$ in the diagonal matrix $\sD$ are given in the unit [1/s] and are related to frequencies $f_{1}$ and $f_{2}$ by $d_1 =  \nicefrac{f_{1}}{Q_1}$ and $d_2 =  \nicefrac{f_{2}}{Q_2}$,
where $Q_1$ and $Q_2$ are the unitless so-called \emph{quality factors}; cf. \cite[(3.38)]{Singiresu1995mechanicalVibration}.
We refer to $d_1$ and $d_2$ also as damping coefficients.
These damping coefficients give the energy decay constant of the harmonic oscillator.
Moreover, we denote with $\lambda>0$ the coupling constant of the system.

%%%%%%%%%%%%%%%%%%%%%%%%%%%%%%%%%%%%%%%%%%
\subsection{Transformation of the system}
\label{sec:Transformation}
%%%%%%%%%%%%%%%%%%%%%%%%%%%%%%%%%%%%%%%%%%

The model in the form of \eqref{eq:original_system} is used in theory to describe a system of coupled oscillators. 
However, it is not the observable system to which we have access in our laboratory experiments. 
In particular, we cannot measure the coupling coefficient $\lambda$ or the frequencies $f_{1}$ and $f_{2}$ directly.
We can only  measure the hybridized modal coordinates $\widetilde{q}_1$ and $\widetilde{q}_2$ of the nanostring that result from a hybridization of $q_1$ and $q_2$.
They do not coincide with the physical coordinates $\bq$ for $\lambda>0$.
To derive the ODEs describing the behavior of the hybridized modal coordinates, we have to perform a coordinate transformation through a matrix $\sT$ from the original physical coordinates to the hybridized ones.
The resulting system of ODEs has a diagonal stiffness matrix and therefore does not explicitly include the coupling coefficient.

To determine the transformation matrix $\sT$, we follow a method similar to the one described in \cite{FrimmerNovotny2014ClassicalBlochEquations}.
Since $\sC$ in \eqref{eq:original_system} is real and symmetric, its eigenvalues are real and we can diagonalize it using an orthogonal matrix \cite[Section 5.4]{Shores2007AppliedLinA}. 
The eigenvalues are given by
\begin{align}
	\eta^2_\pm \coloneqq	
	\frac{1}{2}\left(
	(2\pi f_{1})^2 + (2\pi f_{2})^2 
	\mp \sqrt{4\left(2 \pi \lambda \right)^4 + \left( (2 \pi f_{1})^2-(2\pi f_{2})^2\right)^2}
	\right)
	\label{eq:decoupled_system}
\end{align}
with corresponding orthogonal eigenvectors $v_+=(v^1_+,v^2_+)$ and $v_-=(v^1_-,v^2_-)$. 
Since $\sT$ is orthogonal, it holds that $\sT^{-1}= \sT^\top$. 
We remark that the matrix $\sT$ can be assembled if one knows exactly the values of $f_{1},f_{2}$ and $\lambda$.
However, these values are unknown in the current case and we aim to estimate them by finding the correct transformation $\sT$.

We introduce now the coordinate transformation between the hybridized modal coordinates $\tildebq \coloneqq (\widetilde{q}_1,\widetilde{q}_2)$ and the physical coordinates $\bq$:
\begin{align}
    \tildebq(t) = (\sT \bq)(t).
    \label{eq:variable_transformation}
\end{align}
The system describing the evolution of the hybridized modal coordinates, to which we have access in the experiments, has the following structure:
\begin{subequations}
	\begin{align}
		&\ddot{\tildebq}(t) + \widetilde{\sD}\,\dot{\tildebq}(t) + \widetilde{\sC}\, \tildebq(t) = \widetilde{\bb}(t), 
        \qquad
        t \in [0,\Ttot]
		\label{eq:experimental_model}
		\\
		&\bz = \mathcal{F}(\widetilde\bq|_{[\Ttrans,\Ttot]})
		\label{eq:experimental_outputs}.
	\end{align}
    \label{eq:eigenmodes_model}
\end{subequations}
\noindent
\SecondRevisionComment{The system is completed with the transformed initial condition given by $\tilde{\bq}(0) = \dot{\tilde{\bq}}(0) = 0$.}
The observation $\bz$ is given as the Fourier transform $\mathcal{F}$ of the time-dependent function $\widetilde\bq$ in the time interval $[\Ttrans,\Ttot]$. 
In \eqref{eq:experimental_model}, the matrices containing the damping and eigenfrequency, and the vector containing the external driving are given by
\begin{align}
	\widetilde{\sD} \coloneqq \sT \sD \sT^\top 
    =\begin{pmatrix}
		\widetilde{\sD}_{11} & 	\widetilde{\sD}_{12} \\
		\widetilde{\sD}_{21} & 	\widetilde{\sD}_{22}
	\end{pmatrix}, 
	\qquad
	\widetilde{\sC} \coloneqq  \sT \sC \sT^\top
    = \begin{pmatrix}
		(2\pi\eta_+)^2 & 0 \\
		0 & (2\pi\eta_-)^2
	\end{pmatrix},
    \qquad´
    \widetilde{\bb}(t) = \binom{\widetilde{b}(t)}{\widetilde{b}(t)},
    \label{eq:tildeD_tildeC}
\end{align}
where $\widetilde{b}$ will be specified in \cref{sec:DrivingForces}.
Furthermore, the matrix $\widetilde{\sD}$ has the following structure 
\begin{align}
	\widetilde{\sD} = \sT  \begin{pmatrix}
		2\pi d_1  & 0
		\\
		0 & 2\pi d_2 
	\end{pmatrix} \sT^\top = 
	\frac{2\pi}{v^1_+ + v^1_-} \begin{pmatrix}
		d_1\,v^1_+ + d_2 \,v^1_-\,&\,  (-d_1+d_2)\, v^1_-v^1_+ \\
		(-d_1+d_2)& d_1 \,v^1_- + d_2 \, v^1_+
	\end{pmatrix} .
    \label{eq:Def_tildeD}
 \end{align}
Notice that $\widetilde\sD$ is not diagonal.
The off-diagonal elements can get small if the damping constants $d_1$ and $d_2$ get close to each other.
Furthermore, we see that the diagonal elements in $\widetilde\sD$ are a (weighted) average of the damping constants $d_1$ and $d_2$.

In our physical application, we can measure the eigenvalues $\eta_\pm$ and the quality factors $Q_\pm$ (cf. \cref{sec:FRF_calibration}).
Furthermore, we suppose that the diagonal elements of $\widetilde\sD$ are given by
\begin{align}
    \widetilde\sD_{11} = d_+ = \nicefrac{\eta_+}{Q_+},
    \qquad\qquad
    \widetilde\sD_{22} = d_- = \nicefrac{\eta_-}{Q_-}.
    \label{eq:diagonal_tildeD}
\end{align}

To regain the mathematical model \eqref{eq:original_system} with explicitly including the coupling, we apply $\sT^\top$ from the left to \eqref{eq:experimental_model} and use \eqref{eq:variable_transformation} which results in 
\begin{align}
	\ddot{\bq}(t) + \sT^\top  \widetilde \sD  \sT\,\dot{\bq}(t) + \sT^\top \widetilde\sC \sT \bq(t) = \sT^\top \widetilde{\bb}(t),
    \qquad\qquad
    t \in [0,\Ttot].
	\label{eq:mathematical_model}
\end{align}
For the right-hand side of \eqref{eq:mathematical_model} it holds that
\begin{align*}
	\sT^\top \widetilde{\bb}(t) = \bb(t) 
	= \binom{b_1(t)}{b_2(t)} 
	= \widetilde{b}(t)\binom{v^1_-+v^1_+}{v_-^2+v_+^2}.
\end{align*}

Furthermore, recall that any orthogonal transformation matrix is an element of the orthogonal group $O(2)$ and hence of one of the following forms for $\theta \in (-\pi,\pi)$:
\begin{align}
	\begin{pmatrix}
		\cos(\theta) & -\sin(\theta) \\
		\sin(\theta) & \cos(\theta)
	\end{pmatrix},
	\qquad\qquad
	\begin{pmatrix}
		-\cos(\theta) & \sin(\theta) \\
		\sin(\theta) & \cos(\theta)
	\end{pmatrix}.
	\label{eq:O2_elements}
\end{align}
Thus, we aim to find the correct parameter $\theta \in (-\pi,\pi)$ to calculate the correct transformation $\sT$.
To indicate this fact, we write $\sT_\theta$ for the transformation matrix given the parameter $\theta$.
In the numerical experiments, we use both forms and choose the one that leads to better results.

Notice that we do not have access to the elements in $\sD$ through the measurements.
Hence, we include this in our parameter estimation problem.
In summary, we search for the transformation matrix $\sT_\theta$ and for the correct damping constants $d_1$, $d_2$ in \eqref{eq:original_system}.
From the structure of $\widetilde\sD$, we know that $\widetilde\sD_{ii}$ contains information about $d_i$, $i=1,2$ (cf. \eqref{eq:Def_tildeD}).

The unknowns in this setting are now the transformation parameter $\theta$ for $\sT_\theta$ and the scalars $d_1,d_2$ constituting $\sD$.
To shorten the notation, we introduce the parameter vector containing the parameters we want to estimate as
\begin{align}
    \bp \coloneqq (\theta,d_1,d_2) \in \Rb^3.
    \label{eq:parameter_vector}
\end{align}

We introduce the set of admissible parameters $\sPad$ as
\begin{align}
    \sPad \coloneqq
    \left\lbrace 
    \bp \in \Rb^3 \; \vert \;
    \bp_{\mathsf{min}} \leq \bp \leq \bp_{\mathsf{max}}
    \right \rbrace,
\end{align}
with fixed parameter vectors $\bp_{\mathsf{min}} \leq \bp_{\mathsf{max}}$ (componentwise).
Notice that $\sPad$ is non-empty, convex, and compact in $\Rb^3$.
Incorporating the admissible set $\sPad$ is mainly due to being able to apply standard methods from optimization theory.

%%%%%%%%%%%%%%%%%%%%%%%%%%%%%%%%
\subsection{External driving forces}
\label{sec:DrivingForces}
%%%%%%%%%%%%%%%%%%%%%%%%%%%%%%%

Motivated by \cite{GajoRastelliWeig2019NonlinearCouplingResonator,HuberBelzigWeig2020WeaklyDampedDrivenMode}, we choose the following structure of the time-dependent driving forces (also called external excitation or controls)\textcolor{red}{:}
\begin{align}
	\widetilde{\bb}(t) = \binom{\widetilde{b}(t)}{\widetilde{b}(t)},
	\qquad\qquad
	\widetilde b(t) = A\cos(2\pi u_1 t)+A\cos(2\pi u_2t),
    \qquad\qquad
    t \in [0,\Ttot].
	\label{eq:control_structure}
\end{align}
We consider a control vector $\widetilde{\bb}$ with two identical components. 
This reflects the fact that both modes are exposed to the same excitation signal, assuming that both modes couple equally to the drive.

Here, $A>0$ is the driving amplitude and $u_1,u_2>0$ are \revisionComment{the} driving frequencies.
We choose a fixed driving amplitude $A$ and consider the driving frequencies $u_1$ and $u_2$ as input parameters.

In this setting, we know by trigonometric identities that for our controls it holds that
\begin{align}
	A \cos(2\pi u_1 t) + A \cos(2\pi u_1 t) =
	2 A \cos \left(\frac{2\pi(u_1 + u_2)t}{2}\right) \cos\left(\frac{2\pi(u_1-u_2)t}{2}\right).
\end{align}

\begin{figure}
    \centering
    \begin{subfigure}[l]{0.42\textwidth}
        \includegraphics[width=\textwidth]{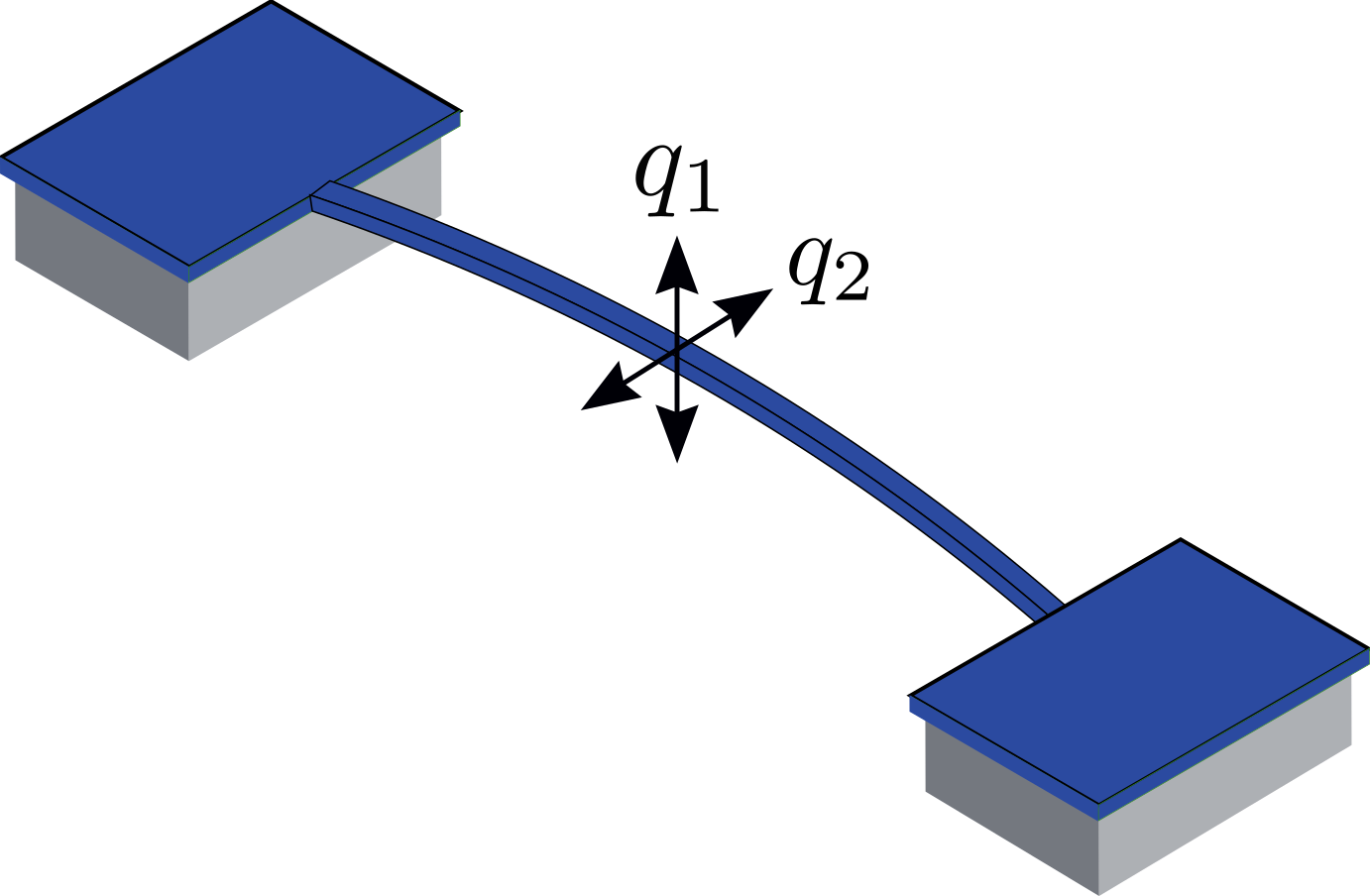}
        \caption{}
    \label{fig:VibrationModes}
    \end{subfigure}
 \hfill
    \begin{subfigure}[l]{0.35	\textwidth}
       \includegraphics[width=\textwidth]{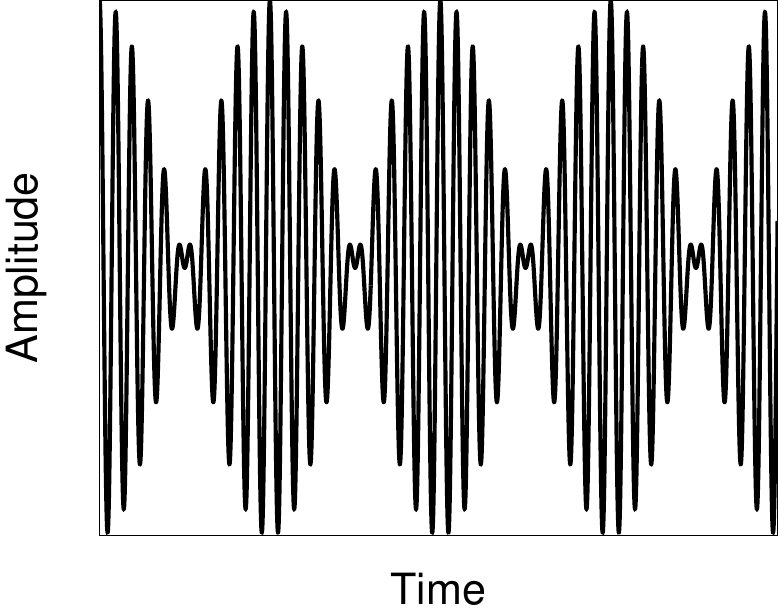}
       \caption{}
       \label{fig:driving}
    \end{subfigure}
    \caption{(a) Vibration modes of the nanostring.
    ;
    (b) Form of external driving (cf. \eqref{eq:external_driving})
    }
    \label{fig:Setting}
\end{figure}

This phenomenon is also known as \emph{beating} \cite[Section 3.3.2]{Singiresu1995mechanicalVibration}; see also \cref{fig:driving}.

Notice that the frequency of one of the Cosine functions on the right-hand side is the difference between the two frequencies $u_1$ and $u_2$.
Since they are close to each other in our case, this results in a slow mode.
Therefore, our time interval in the experiments and simulations should be long enough in order to catch both modes.
\SecondRevisionComment{Moreover, we always choose the pair of driving frequencies in such a way that one is smaller and the other one is larger than the eigenfrequency.}

We choose to have $n_c \in \Nb$ of such driving forces with driving frequency pairs $u_1^m,u_2^m$ for $m=1,\ldots,n_c$. 
To shorten the notation, we introduce the vectors $\bu \in \Rb^2$ representing the input parameters:
\begin{align}
	\bu^m = (u^m_1,u^m_2)^\top,
    \qquad
    m = 1,\ldots,n_c.
    \label{eq:def_epsilon}
\end{align}

Furthermore, we define
\begin{align}
    \widetilde{\bb}^m(t) = (\widetilde b^m(t), \widetilde b^m(t))^\top,
    \qquad\qquad
    \text{ where }
    \qquad
    \widetilde{b}^m(t) =A\cos(2\pi u^m_1 t)+A\cos(2\pi u^m_2t).
    \label{eq:external_driving}
\end{align}

\subsection{Analysis of the steady state solution}

\begin{asmpt}
    \label{asmpt:steady_state}
    The matrix
\[
A_k:=- (2\pi u_k)^2 I + i(2\pi u_k)\widetilde{\sD}+\widetilde{\sC}
\]
is invertible for $k=1,2$ are invertible for all $p \in \sPad$, where $\widetilde{\sD}$ and $\widetilde{\sC}$ are given in \eqref{eq:tildeD_tildeC}.
\end{asmpt}

\begin{lemma}
    \label{lem:parameter-to-state-map}
    \SecondRevisionComment{
   Let \cref{asmpt:steady_state} hold.
	The parameter-to-state map that maps a parameter $\bp \in \sPad$ to the steady-solution of \eqref{eq:eigenmodes_model} is continuously differentiable. }
\end{lemma}

\begin{proof} 
\SecondRevisionComment{
Using $\cos(\omega t)=\tfrac12\left(e^{i\omega t}+e^{-i\omega t}\right)$, 
the forcing can be written as 
$\widetilde{\bb}(t) = \sum_{k=1}^2 \frac{A}{2}\binom{1}{1}e^{i2\pi u_k t} + \frac{A}{2}\binom{1}{1}e^{-i2\pi u_k t}. $
By the linearity of the system and applying the superposition principle, the steady-state solution has the form (see \cite[Chapter 4]{Singiresu1995mechanicalVibration})
\begin{align}
\tilde{\bq}(t) =
\sum_{k=1}^2
\mathbf X_k e^{i2\pi u_k t}
+
\mathbf Y_k e^{-i2\pi u_k t},
\label{eq:steady_state}
\end{align}
with amplitudes $\mathbf X_k$ and $\mathbf Y_k$ for the two modes.
We denote by steady state the long-term behavior of the system after any initial transients have settled.
Substituting \eqref{eq:steady_state} into the differential equation yields
\[
A_k \mathbf X_k=\mathbf F,
\qquad
A_k^\ast \mathbf Y_k=\mathbf F,
\qquad
\mathbf F=\tfrac{A}{2}(1,1)^\top .
\]
Hence
\[
\mathbf X_k=A_k^{-1}\mathbf F .
\]
For $2\times2$ matrices we can directly write the form of the
inverse matrix as $A_k^{-1} = \operatorname{adj}(A_k)(\det A_k)^{-1}$.
The entries of $A_k$ depend polynomially on
$\widetilde{\sD}_{ij}$ and $\eta_\pm$. Hence
$\operatorname{adj}(A_k)$ and $\det A_k$ are polynomial in these parameters.
Therefore
\[
\mathbf X_k
=
\frac{\operatorname{adj}(A_k)\mathbf F}{\det A_k}
\]
is a rational function of the parameters.
Thus $\widetilde{\bq}(t)$ depends smoothly on the entries of
$\widetilde{\sD}$ and on $\eta_\pm$.
}
\end{proof}
\medskip
\SecondRevisionComment{We remark that \cref{asmpt:steady_state} needed for \cref{lem:parameter-to-state-map} is a mild assumption in our context. It is always fulfilled since we assume a non-vanishing damping.
}
\SecondRevisionComment{Moreover, we remark that our numerical algorithm is also able to take arbitrary external drivings into account. Choosing the specific form of \eqref{eq:external_driving} is due to the more precise generation in the laboratory and the analysis of the differentiability of the objective (see \cref{thm:Differentiability_Objective} below).}

%%%%%%%%%%%%%%%%%%%%%%%%%%%%%%%%%%%%%%%%%%%%%%%%%%
\subsection{Drifting eigenfrequencies}
\label{sec:DriftingEigenfrequencies}
%%%%%%%%%%%%%%%%%%%%%%%%%%%%%%%%%%%%%%%%%%%%%%%%%%

An important observation that we made during our work is that the eigenfrequencies $\eta_+$ and $\eta_-$ encounter drifts during the execution of the laboratory experiments, e.g., from variations in the ambient conditions, notably slow temperature fluctuations.
To include this in our model, we assume to know $n_c$ eigenfrequencies that result from linear interpolation between the eigenfrequencies measured at the beginning and the end of our experiments.
For this reason, we define $\widetilde{\sC}^m$ and $\sC^m$ as
\begin{align}
	\widetilde{\sC}^m \coloneqq \begin{pmatrix}
		(2\pi\eta^m_+)^2 & 0 \\
		0 & (2\pi\eta_-^m)^2
	\end{pmatrix},
	\qquad
    \sC^m \coloneqq 
    \begin{pmatrix}
            (2\pi f^m_{1})^2 & -(2\pi \lambda^m)^2
            \\
            -(2\pi \lambda^m)^2 & (2\pi f^m_{2})^2
    \end{pmatrix}
    \qquad\text{for }
    m \in \{1,\ldots,n_c\}.
    \label{eq:matrix_drifted_eigenfrequencies}
\end{align}

We assume that $\sD$ and $\widetilde{\sD}$ do not depend on $m$. 
This is an approximation but we observe in our experiments that this assumption is suitable for the scope of this work.

%%%%%%%%%%%%%%%%%%%%%%%%%%%%%%%%%%%%%%%%%%%%%%%%%%%%%
\section{Physical system and experimental setup}
\label{sec:PhysicalSystem}
%%%%%%%%%%%%%%%%%%%%%%%%%%%%%%%%%%%%%%%%%%%%%%%%%%%%%

\begin{figure}
    \centering
    \begin{subfigure}[l]{0.4\textwidth}
        \includegraphics[width=\textwidth]{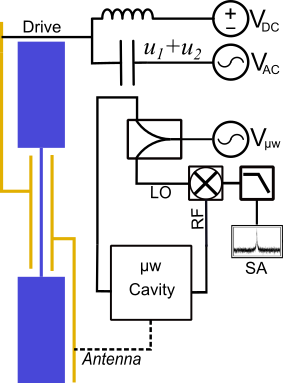}    
        \caption{}
        \label{fig:messsystem}
    \end{subfigure}
    \hfill
    \begin{subfigure}[l]{0.55\textwidth}
       \includegraphics[width=\textwidth]{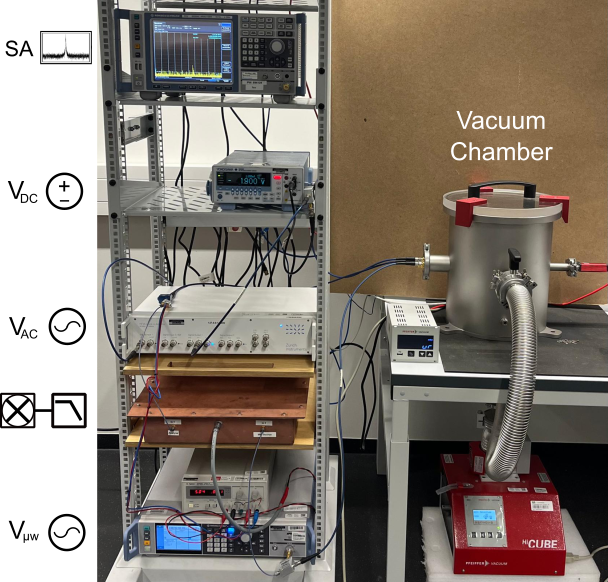}
        \caption{}
        \label{fig:lab}
    \end{subfigure}
    \caption{Experimental measurement setup. (a) Schematic of the measurement setup used to drive and detect the vibration of the nanostring. (b) Photograph of the measurement setup, depicting major electronic components from \cref{fig:messsystem} as well as the vacuum chamber hosting the microwave cavity and the nanostring chip.}
    \label{fig:Setup}
\end{figure}

As described above, the reconstruction method presented in this work aims in particular to identify the coupling parameter $\lambda$ between two orthogonal vibration modes of a nanomechanical string resonator.
The resonator is made of a thin film of amorphous, stoichiometric silicon nitride (Si$_3$N$_4$) deposited on a fused silica wafer, and fabricated using top-down nanofabrication techniques~\cite{BarakatWeig2024CouplingStrength}. 
It has a length of $60\,\mu$m, a thickness of $100$\,nm and a width of $250$\,nm, approximately. 
The experiment is performed at room temperature. 
In order to avoid damping of the nanostring's vibrations from the surrounding medium, it is placed in a vacuum chamber at a pressure below $10^{-4}$\,mbar (see \cref{fig:lab}).

The vibrational modes of the nanostring exhibit two orthogonal \enquote{polarizations} which are associated with a vibration along the out-of-plane and in-plane direction; cf. $q_1$ and $q_2$ in \cref{fig:VibrationModes}. 
In principle, the nanostring, being a continuous system, exhibits an infinite number of vibrational modes in each of these polarization direction. 
In the context of this work, we will only be interested in its two fundamental modes. This is justified as the two fundamental modes exhibit similar eigenfrequencies which are spectrally well separated from all other higher-order harmonics.
The lowest-lying out-of-plane and in-plane flexural mode $q_1$ and $q_2$, schematically depicted in \cref{fig:VibrationModes}, are hybridized into $\tildeq_1$ and $\tildeq_2$ by the coupling $\lambda$, where these vibration coordinates are rotated, but also orthogonal to each other~\cite{meirovitch2010fundamentals}. 
The model describing the evolution in time of the two vibrational modes is given in \eqref{eq:eigenmodes_model}.

\begin{figure}
	\begin{subfigure}[l]{0.45\textwidth}
		\includegraphics[width=\textwidth]{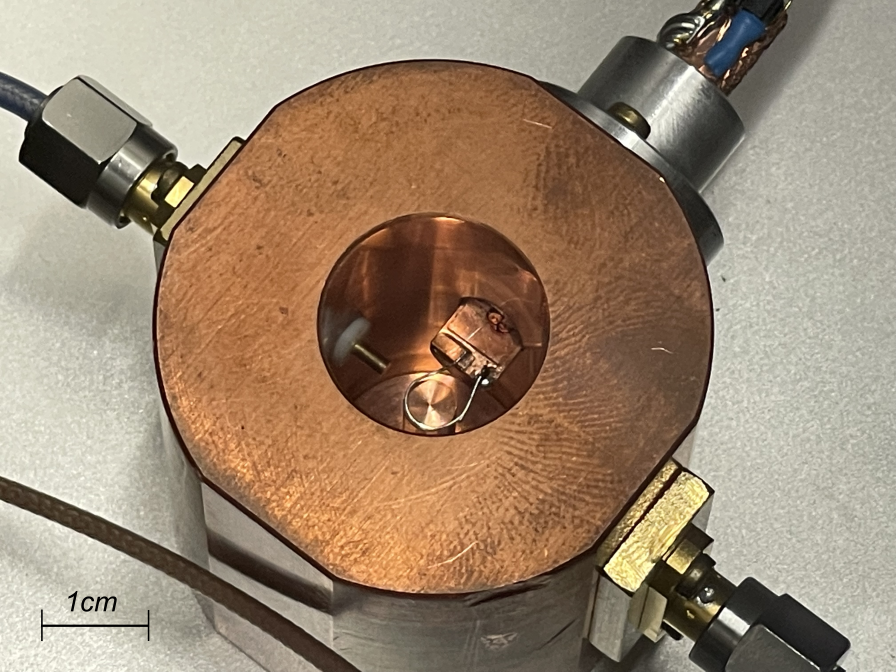}
		\caption{}
		\label{fig:cavity}
	\end{subfigure}
\hfill
	\begin{subfigure}[l]{0.5\textwidth}
		\includegraphics[width=\textwidth]{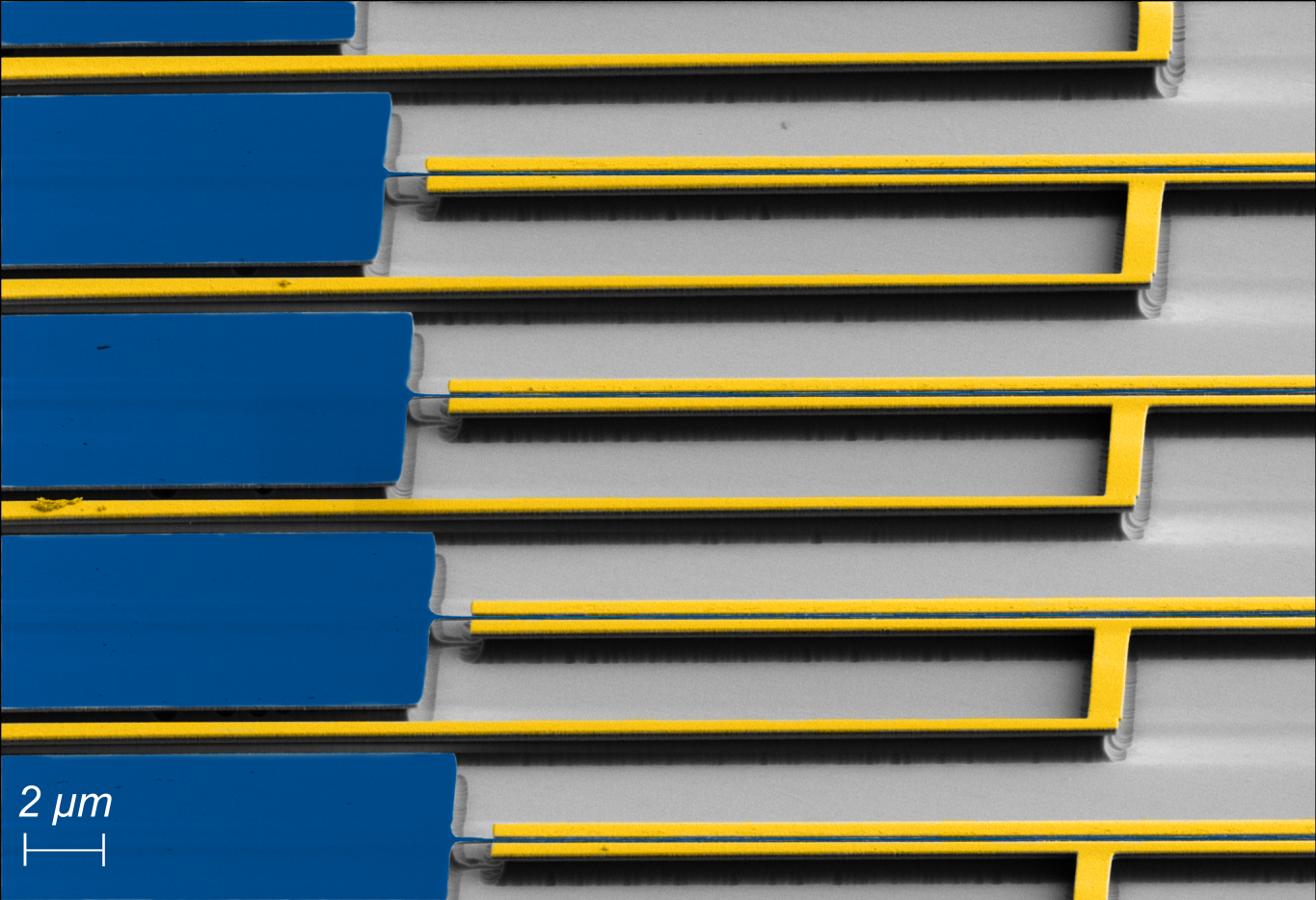}
		\caption{}
		\label{fig:string}
	\end{subfigure}
	\caption{Microwave cavity and nanostring chip. 
    (a) Photograph of the coaxial quarter-wave microwave cavity. 
    The loop antenna coupling the cavity mode to the electrical field surrounding the nanostring as well as the (transparent) nanostring chip can be discerned above and next to the center stub, respectively. 
    (b) Scanning electron micrograph of a representative device depicting a partial view of several nanostrings (blue) between pairs of electrodes (yellow).}
	\label{fig:Laboratory}
\end{figure}

Let us now describe the methods used for the vibrational excitation and measurement of the nanostring, which are illustrated in \cref{fig:Setup} and \cref{fig:Laboratory}.
The vibration is excited by means of an electric field surrounding the nanostring, induced by voltages applied between two electrodes positioned on either side as shown in \cref{fig:messsystem} and \cref{fig:string}. 
Since the nanostring is made of a dielectric material, the surrounding electric field will cause an electric polarization within, forming a collective electric dipole, which then follows the known dynamics of a dipole in an inhomogeneous electric field. 
The inhomogeneous electric field is engineered by placing the string out of the field's axes of symmetry.
In this way, a net dipole force is created in both vibrational coordinates, which then allows to induce the mechanical motion by combining a static (DC) with a dynamic (AC) voltage as depicted in \cref{fig:Setup}; see \cite{unterreithmeier2009universal}. 
In this experiment, the dynamic voltage $V_{\mathrm{AC}}$ consists of the sum of two excitation signals with two frequencies $u_1$ and $u_2$ close to each other (within the mode's linewidth) according to \eqref{eq:control_structure}.
Notice that $V_{\mathrm{DC}}$ further determines the eigenfrequencies as well as other parameters of the system such as the linear coupling between the two modes and their damping rates, and is chosen according to the 
desired coupling strength.

The measurement of the string's vibrations relies on the detection of the modulation in the capacitance between the two electrodes, which is caused by the vibration of the dielectric string in the inhomogeneous electric field.
To resolve the resulting minute capacitance modulation, a microwave-cavity assisted heterodyne detection scheme is employed (see~\cite{le_3d_2023} for more detail).
We employ a three-dimensional quarter-wave coaxial cavity that couples to the electric field between the electrodes via a loop antenna (see \cref{fig:messsystem} and \cref{fig:cavity}). 
The cavity output signal (RF) is mixed with a reference signal (LO), amplified and filtered, before being recorded using a spectrum analyzer (see \cref{fig:Setup}).

A salient feature of our experimental setup is the tunability of the system with $V_{\mathrm{DC}}$. In this work, we fixed the system's parameters by selecting a specific bias voltage, $V_{\mathrm{DC}}=-12$\,V. 
At this operating point, the two modes exhibit a sizable coupling strength~\cite{BarakatWeig2024CouplingStrength}. 
In addition, their eigenfrequencies have been tuned to a situation where the transformation matrix $\sT$ is clearly distinct from the identity matrix, such that the coupling takes effect on the vibrational dynamics. 
It is in fact one advantage of our strategy that we can estimate the coupling parameter $\lambda$ for a single $V_{\mathrm{DC}}$ voltage without the need to acquire information about the system at other $V_{\mathrm{DC}}$ voltages; see also \cref{sec:Conclusion,sec:Comparison}.

%%%%%%%%%%%%%%%%%%%%%%%%%%%%%%%%%%%%%%%%%%%%%%%%%%%
\section{Scaling of simulation and experiment}
\label{sec:FRF_calibration}
%%%%%%%%%%%%%%%%%%%%%%%%%%%%%%%%%%%%%%%%%%%%%%%%%%%

To ensure that the observations of the simulation and experimental setup are comparable, we need to scale the output of the simulations. 
The reason why we need this is that the physical quantities in the laboratory experiments are all measured in the unit [V] (volt).
To determine the scaling factor, we measure the \emph{frequency response function} (FRF) in the laboratory and compare it with the result of the computer simulations. 
Experimentally, the FRF for a specific mode is obtained by exciting the nanostring with several frequencies close to the eigenfrequency of the mode.
At each excitation frequency, the amplitude response of the nanostring is recorded.
This procedure is also known as \emph{sweeping}.
Then we use a software to determine the resonance frequency $\eta_\pm$, the amplitude $A_\pm$ of the peak, and the quality factor $Q_\pm$.
To compute the scaling, the shape of the squared FRF is important. 
This should have a \emph{Lorentzian} shape with a maximum peak $A_\pm$ at the resonant frequency $\eta_\pm$; this shape is also called \emph{Cauchy distribution} in the realm of probability density functions; see, e.g., \cite{johnson1972continuousUnivariateDistributions}.
The scaling factor is now the ratio between the peak amplitudes of the displacements in equilibrium in the laboratory experiment and the simulations.

This Lorentzian shape is defined by the following formula
\begin{align}
	\rho_{L}(f;\eta_\pm,d_\pm,A_\pm) 
	\coloneqq
	\frac{A_\pm - \xi}{1+ \left(\frac{2\pi f-2\pi \eta_\pm}{d_\pm}\right)^2}
	\qquad\text{for }f >0,
	\label{eq:Lorentzian-shape}
\end{align}
where $d_\pm$ is a scale parameter that specifies the full-width at half-maximum (FWHM) and equals the damping parameter.
As introduced in \eqref{eq:diagonal_tildeD}, we have the relation $Q_\pm=\nicefrac{\eta_\pm}{d_\pm}$.
Furthermore, $0<\xi \ll 1$ is the noise level of the lab data.

\begin{remark}
    An important observation is that the scaling also depends on the unknown damping coefficients $d_1,d_2$ and transformation matrix $\sT_\theta$; see also \cref{algo:Recon} below.\hfill$\Diamond$
\end{remark}

We summarize our procedure to get our $\bp$-dependent scaling factor $\chi_\bp$ in \cref{algo:FRF_Calibration}.
\SecondRevisionComment{Notice that we do not need Fourier transforms here since we stay the time-domain.}

\bigskip

\SecondRevisionComment{
\begin{lemma}
    \label{lem:FRF_smooth}
    Let \cref{asmpt:steady_state} hold.
    Let $\widetilde{\bq}$ be the solution of \eqref{eq:experimental_model} using $\bp \in \sPad$ and $\cos(2\pi f_S^j t)$ as input.
    Then the amplitude of $\widetilde q_1$, given as $\max_{t \in \SecondRevisionComment{[\Ttrans,\Ttrans+2\pi\,f_S^j]}} |\widetilde q_1(t;\bp)|$, depends smoothly on $\bp$.
\end{lemma}
}

\begin{proof}
	\SecondRevisionComment{We argue similar to the proof of \cref{lem:parameter-to-state-map} since the driving force here is also a cosine with a given frequency. Hence, the maximum of the amplitude depends smoothly on the parameters.}
\end{proof}

\begin{algorithm}
	\caption{(FRF Calibration Laboratory/Simulation)}
	\label{algo:FRF_Calibration}
		\begin{algorithmic}[1]
            \Require parameter vector $\bp \in \Rb^3$, 
            FRF data set from laboratory experiment $\boldsymbol C_L \in \Rb^{N_{\mathrm{FRF,L}}}$, frequency discretization $\{f_S^1,\ldots,f_S^{N_{\mathrm{FRF,S}}}\}$ for simulation of sweeping
            \State Initialize empty list $\boldsymbol C_S \gets \{\}$;
            \For{$j=1,\ldots,N_{\mathrm{FRF,S}}$}
                \State Simulate the system \eqref{eq:experimental_model} using $\bp$ and $\cos(2\pi f_S^j t)$ as input, 
                \SecondRevisionComment{denote the solution for the first mode by $\tilde{q}_1(t;\bp)$}
                \State Append the maximum amplitude $\max_{t \in \SecondRevisionComment{[\Ttrans,\Ttrans+2\pi\,f_S^j]}} |\widetilde q_1(t;\bp)|$ of $\widetilde q_1$ to $\boldsymbol C_S$;
            \EndFor
            \State Set $\chi_\bp \gets \nicefrac{\max_{1 \leq i \leq N_{\mathrm{FRF,L}}} (\boldsymbol C_L)_i}{\max_{1 \leq i \leq N_{\mathrm{FRF,S}}}(\boldsymbol C_S)_i}$;
			\\
            \Return $\chi_\bp$
		\end{algorithmic}
\end{algorithm}

%%%%%%%%%%%%%%%%%%%%%%%%%%%%%%%%%%%%%%%%%%%%%%%
\section{Parameter identification problem}
\label{sec:ParamIdentProblem}
%%%%%%%%%%%%%%%%%%%%%%%%%%%%%%%%%%%%%%%%%%%%%%%

Recall that we aim at estimating the damping coefficients $d_1,d_2$ in $\sD$ and transformation matrix $\sT_\theta$. 
The latter one will lead to estimates of the matrices $\sC^m=\sT_\theta^\top\,\widetilde\sC^m\,\sT_\theta$ for $m=1,\ldots,n_c$, where $\widetilde{\sC}^m$ is known (cf. \cref{sec:Setting}).
To stress the dependence of $\sD$ and $\sT_\theta$ on the unknown parameter vector
\begin{align*}
    \bp=(\theta,d_1,d_2)
\end{align*}
defined in \eqref{eq:parameter_vector}, we use the notation $\sT_\bp \coloneqq \sT_\theta$ and
\begin{align*}
    \sD_\bp \coloneqq \begin{pmatrix}
        2\pi d_1 & 0 \\
        0 & 2\pi d_2
    \end{pmatrix},
    \qquad\qquad
    \widetilde{\sD}_\bp \coloneqq \sT_\bp \sD_\bp \sT_\bp^\top.
\end{align*}
We have the following structure of the ODEs in the final identification:
\begin{align}
    \ddot{\bq}^m(t) + \sD_\bp \, \dot{\bq}^m(t) + \sT_\bp^\top \, \widetilde{\sC}^m \, \sT_\bp\, \bq^m(t) = \sT_\bp^\top \widetilde{\bb}^m(t),
    \quad\qquad
    t \in [0,\Ttot],
    \qquad\qquad
    m \in \{1,\ldots,n_c\}
	\label{eq:ODE_final_ident}
\end{align}
with known $\widetilde{\sC}^m, \widetilde{\bb}^m$ and unknown $\sT_\bp$, $\sD_\bp$. 
\SecondRevisionComment{Since the The initial conditions are not changed by this transformation, we have $\bq(0)=\dot{\bq}(0) = 0$.}

Notice that in \eqref{eq:ODE_final_ident} the untransformed matrix $\sD_\bp$ and the transformed matrix $\widetilde{\sC}^m$
appear. 
The reason for this is given by the limitations and possibilities of the measurements. 
More specifically, we can measure all entries of the diagonal matrix $\widetilde{\sC}^m$ but only the diagonal elements of $\widetilde\sD_\bp$. 
However, $\widetilde\sD_\bp$ is not diagonal (cf. \eqref{eq:Def_tildeD}).
Hence, we reconstruct the untransformed matrix $\sD_\bp$.

Next, we denote by $\bz_\bp^{\bu^m} \coloneqq \mathcal{F}(\sT_\bp \, \bq^m|_{[\Ttrans,\Ttot]})$, $m=1,\ldots,n_c$, the Fourier transformations within the time interval $[\Ttrans,\Ttot]$ of the solutions of $\eqref{eq:ODE_final_ident}$ with transformation $\sT_\bp$ of the structure \eqref{eq:O2_elements} generated with $\theta$, using $d_1,d_2$ as damping coefficients, and applying the control parameters $\{\bu^m\}_{m=1}^{n_c}$ (cf. \eqref{eq:def_epsilon}).
Moreover, we define $\{\bz_\star^{\bu^m}\}_{m=1}^{n_c}$ as the experimental data generated in the laboratory using the same control parameters $\{\bu^m\}_{m=1}^{n_c}$.

Our goal is to minimize the difference between the two peaks of the Fourier transformation.
The reason that there appear exactly two peaks is that we are in the linear setting. 
Here, our solution is the superposition of harmonic oscillators under a periodic force for which the behavior is well-known and given as the sum of periodic functions; see, e.g., \cite[Chapter 4]{Singiresu1995mechanicalVibration}.
\revisionComment{Hence, within our linear setting, we cannot change the locations of the two peaks in equilibrium since they are fixed by the control/driving frequency.
However, by optimizing our parameter vector $\bp$, we can adjust the amplitude of the peak such that simulation and experiment coincide better.}

Furthermore, we need to scale the simulation in order to be able to compare the amplitudes of the simulations with the ones of the laboratory experiments; see \cref{sec:FRF_calibration} above.
The procedure to get the scaling factor $\chi_{\bp}$ is described in \cref{algo:FRF_Calibration}.
Notice that $\chi_\bp$ depends on $\bp$. 
Therefore, we need to execute \cref{algo:FRF_Calibration} each time the coefficient vector $\bp$ is changed.

We consider the following definition of the deviation between simulation and experiment.

\begin{definition}
	\label{def:deviation}
	Let $z_\bp^{u}$ \revisionComment{$\in \mathbb{R}^+$} be the amplitude \revisionComment{of the absolute value of the} Fourier transform of the simulations using the parameter vector $\bp$ and the driving frequency $u>0$.
	Furthermore, let $z_\star^u$ be the corresponding amplitude in the experiments.
	Then, we define the \emph{relative deviation} as
	\begin{align}
		\mathfrak{e}(\, z_\bp^{u},z_\star^u) 
		\coloneqq \frac{|z_\bp^{u} - z_\star^u|}{z_\star^u} ,
		\label{eq:deviation_def}
	\end{align}
    where $|\cdot|$ denotes the absolute value.
	Furthermore, we define the vector containing all relative deviations
	\begin{align}
		\mathfrak{E}_{n_c}(\bz_\bp^{\bu},\bz_\star^{\bu}) \coloneqq 
        \left(
        \begin{array}{c}
            \mathfrak{e}(z_\bp^{u_1^m}, z_\star^{u_1^m})\\[1mm]
            \mathfrak{e}(z_\bp^{u_2^m}, z_\star^{u_2^m})
        \end{array}
        \right)_{m=1}^{n_c}\in\Rb^{2n_c}.
		\label{eq:deviation_vector}
	\end{align}
\end{definition}

This definition of the deviation is used in the objective $\sJ$ below.
\revisionComment{Numerically, we get $z^u_\bp$ by taking the value at the driving frequencies since we know that the maximum is located there; see also \cref{sec:Solving_ODEs}.}

Furthermore, we choose some references for the parameters in the vicinity of where we expect the optimal solution to be located.
In particular, we take $\bp_1^\rmref = \theta^\rmref=\pi/2+\pi/8$ since we expect from experimental results that the correct coupling has the corresponding order of magnitude. 
Moreover, we choose $\bp_2^\rmref = d_1^\rmref=\nicefrac{\widetilde\sD_{11}}{2\pi}$, and $\bp_3^\rmref = d_2^\rmref=\nicefrac{\widetilde\sD_{22}}{2\pi}$, since we do not expect too much change in the damping constants compared to measurements $\widetilde{\sD}_{11}$ and $\widetilde{\sD}_{22}$.

Our objective is given  by the sum of the two following functions. 
\begin{align}
    \sJ_{\mathrm{fit}}(\bp,\bz_\bp,\chi_\bp) 
    \coloneqq \sum_{m=1}^{n_c} \mathfrak{e}(\chi_\bp \, z_\bp^{u_1^m},z^{u_1^m}_\star) 
    + \mathfrak{e}(\chi_\bp \, z_\bp^{u_2^m},z^{u_2^m}_\star) ,
\qquad
    \sJ^\nu_{\mathrm{reg}}(\bp) \coloneqq  \frac{\nu}{2} 
        \,{\|\bp - \bp^\rmref\|}_2^2.
        \label{eq:summand_objective}
\end{align}

\SecondRevisionComment{Notice that $\sJ_{\mathrm{fit}}$ is a function of the parameters, the steady state solution and the scaling factor. Moreover, the steady state and the scaling factor depend (smoothly) on the parameters. We introduce the reduced functional $\sJ$, only depending on the parameters $\bp \in \sPad$}.
Now the identification (or parameter estimation) problem is given by
\begin{align}
    %\tag{\textbf P}
    \min 
    \sJ(\bp)
    \coloneqq
    &\sJ_{\mathrm{fit}}(\bp,\bz_\bp,\chi_\bp)
    \,+\,  \sJ^\nu_{\mathrm{reg}}(\bp)
    \quad
    \text{subject to}\quad\bp \in \sPad, 
    \label{eq:online_step_minimization}%
\end{align}
where $\bz_\bp  = \mathcal{F}(\sT_\bp \, \bq_\bp|_{[\Ttrans,\Ttot]})$, $\bq_\bp$ fulfills \eqref{eq:ODE_final_ident} and $\chi_\bp$ is the scaling factor (cf. \cref{sec:FRF_calibration}) for $\bp=(\theta,d_1,d_2)$. \SecondRevisionComment{Notice that $\chi_\bp$ depends on $\bp$ since in the FRF the coupling is inherent (see \cref{sec:FRF_calibration}). This is a crucial property of the optimization problem.}
Problem \eqref{eq:online_step_minimization} is a nonlinear, nonconvex inverse problem with a Tikhonov regularization parameter $\nu\geq 0$ \cite{HaberAscherOldenburg2000OptimizationNonlinInverse}.

We remark that there exists a solution to \eqref{eq:online_step_minimization} by standard arguments; see, e.g., \cite{BartschDenkVolkwein2023adjointCalibrationSDE,BurgerMuhlhuber2002IterativeRegularizationParameterIdentification,Kirsch2021IntroInverseProblems,KunischZou1998IterativeRegularizationInverseProblems}.

\subsection{\SecondRevisionComment{Analysis of the objective functional}}

\SecondRevisionComment{In this section, we discuss the regularity of $\sJ$.}

\begin{theorem}
	\label{thm:Differentiability_Objective}
    \SecondRevisionComment{The functional $\sJ:\sPad \rightarrow \mathbb{R}$ defined in \eqref{eq:online_step_minimization} is continuously differentiable with respect to $\bp$ except at points where $z_\bp^{u} \neq z_\star^u$. }
\end{theorem}

\begin{proof}
\revisionComment{
The functional consists of two parts. 
The part $\mathscr{J}_{\mathsf{reg}}^\nu$ is continuously differentiable as a composition of continuously differentiable functions.
The part $\mathscr{J}_{\mathsf{fit}}$ is more delicate and involves a composition of many functions and operators.\\
Notice that the map from the parameter $\boldsymbol{p}$ to the solution $\bq$ of \eqref{eq:ODE_final_ident} is linear and therefore smooth.
Furthermore, the application of the Fourier transform is linear and hence smooth.}
\\
\SecondRevisionComment{
The amplitude of the solution depends explicitly on the parameters and is smooth in our setting by means of \cref{lem:parameter-to-state-map}. Hence, the map  $\bp \rightarrow \bz_\bp$ is smooth.
Furthermore, by virtue of \cref{lem:FRF_smooth}, the map $\bp \rightarrow \chi_\bp$ is smooth.
Hence, $\chi_\bp \bz_\bp$ is smooth as a procudct of smooth functions. \\
It remains to discuss the regularity of the deviation function $\mathfrak{e}$ defined in \cref{def:deviation}.
Taking the (relative) difference in absolute value of the maximum amplitudes calculated from the simulation and those measured in experiment is in general only Lipschitz continuous. 
Hence, by Rademacher's theorem \cite[Section 3]{EvansGariepy2015MeausureTheory}, it is differentiable almost everywhere. 
However, by the properties of the absolute value function, it is smooth outside the origin. In the current setting, this means that it is smooth if $z_\bp^{u} \neq z_\star^u$.} 
\end{proof}

\medskip

\SecondRevisionComment{We remark that the points $z_\bp^{u} = z_\star^u$ are the best possible solutions of the unregularized problem  $\nu \rightarrow 0$.
More specifically, the functional $\sJ$ is bounded from below by zero.
For $\nu=0$, the only points for a vanishing functional are $z_\bp^{u} = z_\star^u$.
Hence, the functional is only nonsmooth if its lower bound is reached.
}

\subsection{\SecondRevisionComment{Optimisation Strategy}}
To solve \eqref{eq:online_step_minimization}, we apply an iterative strategy commonly used to tackle inverse problems.
First, we choose an initial guess $\bp^0 \in \sPad$ for the parameter vector and an initial Tikhonov regularization parameter $\nu_0$.
Then, we solve \eqref{eq:online_step_minimization} with these quantities, shrink the regularization parameter, and check a termination criterion.
If it is not fulfilled, we start again solving \eqref{eq:online_step_minimization} with the updated parameter vector as the initial guess and with the updated regularization parameter.
\SecondRevisionComment{Notice that every time we update $\bp$, we also have to calculate $\chi_\bp$ again using \cref{algo:FRF_Calibration}.}
We summarize the procedure in \cref{algo:Recon}.

\begin{algorithm}
	\caption{(Iterative reconstruction procedure)}
	\label{algo:Recon}
		\begin{algorithmic}[1]
			\Require Laboratory data $\{\bz^{\bu^m}_\star\}_{m=1}^{n_c}$, $\ell_{\max} \in \Nb$, tolerance $\text{tol}>0$, regularization parameter update $\beta \in (0,1)$;
			\State Initial guess $\bp^0 \in \sPad$, $\ell \gets 0$, initialize $E \gg \text{tol}$, initialize regularization parameter $\nu_0 \gg 0$;
			\While {$E>\text{tol}$ \textbf{and} $\ell < \ell_{\max}$}
				\State Solve \eqref{eq:online_step_minimization} %with \texttt{fmincon} 
                \SecondRevisionComment{with an interior-point method}
                using initial guess $\bp^\ell$ and $\nu=\nu_\ell$  to get $\bp^{\ell+1}$ \SecondRevisionComment{ which requires running \cref{algo:FRF_Calibration} in each iteration}
                \State Set $\nu_{\ell+1} = \beta \, \nu_\ell$ \Comment{update regularization parameter}
				\State Set $E \gets \|\bp^{\ell+1} - \bp^\ell \|_2  +|\sJ_{\mathrm{fit}}-\sJ_{\mathrm{reg}}^{\nu^{\ell+1}}|$;
                \State Set $\ell \gets \ell+1$ 
			\EndWhile
			\\ \Return $\bp^\ell$
		\end{algorithmic}
\end{algorithm}

\begin{figure}
	\centering	
	\begin{tikzpicture}[x=0.75pt,y=0.75pt,yscale=-0.95,xscale=0.95]
		\node[draw=none,fill=none] at (400,100){\includegraphics[scale=0.8]{figs/Setup/Resonator.png}};
		
		%Shape: Rectangle [id:dp9553497766596428] 
		\draw   (17,30) -- (114,30) -- (114,71) -- (17,71) -- cycle ;
		%Shape: Rectangle [id:dp6910585410959749] 
		\draw   (118,230) -- (190,230) -- (190,290) -- (118,290) -- cycle ;
		%Shape: Rectangle [id:dp8278566941376] 
		\draw   (255,231) -- (381,231) -- (381,291) -- (255,291) -- cycle ;
		%Shape: Rectangle [id:dp8555360366516629] 
		\draw   (448,231) -- (642,231) -- (642,290) -- (448,290) -- cycle ;
		%Right Arrow [id:dp010422443156400618] 
		\draw   (128,44) -- (170,44) -- (170,39) -- (198,50) -- (170,61) -- (170,55) -- (128,55) -- cycle ;
		
		%Straight Lines [id:da06399749261734677] 
		\draw    (191.42,260.25) -- (252.42,260.73) ;
		\draw [shift={(254.42,260.75)}, rotate = 180] [color=black][line width=0.75]    (10.93,-3.29) .. controls (6.95,-1.4) and (3.31,-0.3) .. (0,0) .. controls (3.31,0.3) and (6.95,1.4) .. (10.93,3.29)   ;
		
		%Straight Lines [id:da7507397841764643] 
		\draw    (381.42,260.75) -- (444,260.99) ;
		\draw [shift={(446,261)}, rotate = 180] [color=black][line width=0.75]    (10.93,-3.29) .. controls (6.95,-1.4) and (3.31,-0.3) .. (0,0) .. controls (3.31,0.3) and (6.95,1.4) .. (10.93,3.29)   ;
		
		%Right Arrow [id:dp7807426419827779] 
		\draw   (324.5,156) -- (324.5,198) -- (330,198) -- (319,226) -- (308,198) -- (313.5,198) -- (313.5,156) -- cycle ;
		%Bend Up Arrow [id:dp8788058136776817] 
		\draw   (114,261.25) -- (63.5,261.25) -- (63.5,88) -- (75,88) -- (63.5,71) -- (52,88) -- (63.5,88) -- (63.5,261.25) -- (114,261.25) -- cycle ;
		
		% Text Node
		\draw (31,38) node [anchor=north west][inner sep=0.75pt]   [align=center] {Drive Signals\\ \eqref{eq:control_structure}};
		% Text Node
		\draw (135,247) node [anchor=north west][inner sep=0.75pt]   [align=center] {Model\\\eqref{eq:ODE_final_ident}};
		% Text Node
		\draw (276,244) node [anchor=north west][inner sep=0.75pt]   [align=center] {\\Reconstruction
			\\\cref{algo:Recon}};
		% Text Node
		\draw (478,244) node [anchor=north west][inner sep=0.75pt]   [align=center] {\\Solution of parameter\\estimation problem \eqref{eq:online_step_minimization}};
		% Text Node
		\draw (132,70) node [anchor=north west][inner sep=0.75pt]   [align=center] {$\displaystyle u^m_1, u^m_2$};		
		
		\node[draw=none,fill=none] at (160,1){\includegraphics[scale=0.15]{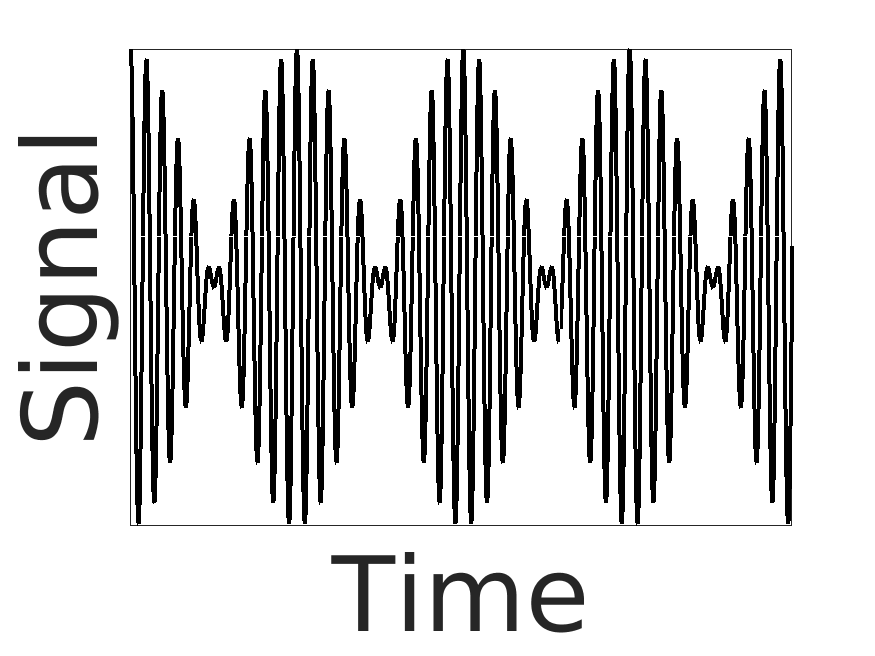}};
		\node[draw=none,fill=none] at (250,190){\includegraphics[scale=0.15]{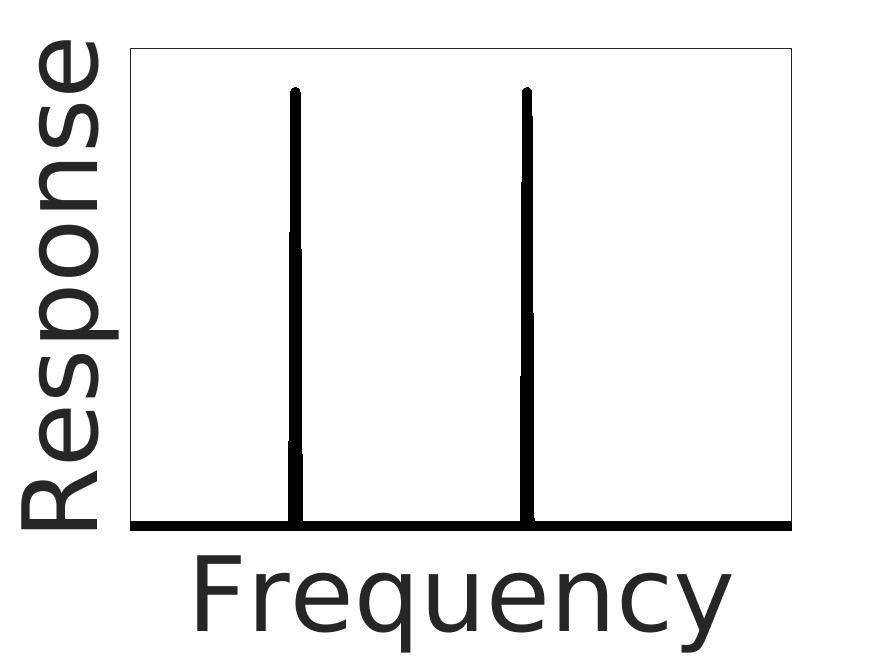}};
	\end{tikzpicture}
	\caption{Schematic visualization of the proposed strategy.}
	\label{fig:Schema}
\end{figure}

A schematic visualization of our whole strategy to estimate the coupling and damping coefficients in a system of oscillators is given in \cref{fig:Schema}.
Starting from our model \eqref{eq:ODE_final_ident}, we apply the same drive signals in the experiment and in simulations.
Finally, the results are used in \cref{algo:Recon} in order to solve the parameter estimation problem \eqref{eq:online_step_minimization}.

%%%%%%%%%%%%%%%%%%%%%%%%%%%%%%%%%%%%%%%%%%%%%%%
\section{Numerical implementation}
\label{sec:Solving_ODEs}
%%%%%%%%%%%%%%%%%%%%%%%%%%%%%%%%%%%%%%%%%%%%%%%%

In this section, we describe our strategy to numerically solve \eqref{eq:ODE_final_ident} over the long time horizon $\Ttot$.
This is needed for two reasons: 
First, we have to integrate over the transient phase $\Ttrans$ since we need to reach a steady state; 
second, we need a sufficiently long signal to detect the correct frequencies numerically in the Fourier transformation; 
third, we have to deal with the beating phenomenon described in \cref{sec:DrivingForces} and hence need a sufficiently large $T$  to correctly resolve all frequencies in the signal.

Since we consider a system of linear ODEs, we can apply the \emph{Laplace transformation} $\mathcal{L}$ in order to get its analytical solution. 
For fixed $m \in \{1,\ldots,n_c\}$, the Laplace transformation $\widehat{q}^m_i(s)=\mathcal{L}\{{q}_i^m\}(s)$ ($i=1,2$) is given as the solution of
\begin{align}
	\left(
	\begin{array}{cc}
		s^2 +  (\sD_\bp)_{11} \, s + (\sC^m_\bp)_{11} & (\sD_\bp)_{12} \, s + (\sC^m_\bp)_{12} \\
		(\sD_\bp)_{21} \, s + (\sC^m_\bp)_{21} & s^2 + (\sD_\bp)_{22} \, s + (\sC^m_\bp)_{22} \\
	\end{array}
	\right)
	\binom{\widehat{q}_1(s)}{\widehat{q}_2(s)}
 =
	\left(
	\begin{array}{c}
		\big((\sT_\bp)_{11} + (\sT_\bp)_{21} \big) \, \widehat{b}^m(s) \\
		\big((\sT_\bp)_{12}+(\sT_\bp)_{22} \big) \, \widehat{b}^m(s) \\
	\end{array}
	\right),
	\label{eq:transformations_problem}
\end{align}
with coordinate transformation matrix $\sT_\bp$, $\sC^m_\bp\coloneqq\sT_\bp^\top \widetilde{\sC}^m\sT_\bp$ and 
\begin{align*}
    \widehat{b}^m(s) 
    \coloneqq \mathcal{L}\{\widetilde{b}^m\}(s) = A \, \left(\frac{ 2\pi u^m_1}{(2\pi u^m_1)^2+s^2}
    +\frac{2\pi u^m_2}{(2\pi u^m_2)^2+s^2}\right).
\end{align*}
After solving \eqref{eq:transformations_problem}, we utilize the built-in Matlab function \texttt{ilaplace} to compute $\bq$ and the function  \texttt{rectangularPulse} in the time interval $[T_{\text{trans}},T_{\text{trans}}+T]$ to have a signal of finite length and discarding the transient phase.

The last step is to apply the built-in Matlab function \texttt{fourier} in order to get the Fourier transformation of the function. 
We then evaluate the Fourier transformation at the frequencies at which we also have experimental data.
For the generation of the frequency interval in which we want to evaluate the Fourier transformation, we split the time interval $[T_{\text{trans}},T_{\text{trans}}+T]$ into $N_t$ equidistant cells of length $\Delta t = \nicefrac{T}{N_t}$.
The frequency interval is then given by
$[f_{\rm a}, f_{\rm b}] = \frac{1}{\Delta t \, T}\left[0, \nicefrac{T}{2}\right]$, where $T$ is the length of the time interval.
In \cref{tab:simulation_hyper_parameters_damping_coupling}, we present the parameters that we use for the simulation.
Moreover, we use $\bp_{\min}=[-2\pi,-0.1,-0.1]$ and $\bp_{\max}=[2\pi,0.1,0.1]$.
\revisionComment{The extraction of the two peaks in the Fourier transform can be performed without solving another optimization problem since we know the position of the amplitudes given by the driving frequencies.}

\begin{table}
    \centering
    \begin{tabular}{c|c||c|c||c|c}\toprule
        $\Delta t$ & $10^{-6}$ &
        $T$  & $3\cdot10^4$ &
        $T_{\mathrm{trans}}$ & $2\cdot10^3$ 
        \\
        \midrule
        $\beta$ & $0.1$ & 
        $\nu_0$ & $10^{-1}$ & 
       $ $ &  \\\bottomrule
    \end{tabular}
    \caption{Simulation parameters and hyperparameters. }
    %https://gitlab.inf.uni-konstanz.de/jan.bartsch/oscillators/-/commit/c01f02a0412934700c04da476a5f086c99030a62
    \label{tab:simulation_hyper_parameters_damping_coupling}.
\end{table}

%%%%%%%%%%%%%%%%%%
\section{Identification using the proposed framework}
\label{sec:NumExp}
%%%%%%%%%%%%%%%%%%%

In this section, we present the results of our proposed framework obtained with laboratory data to validate its effectiveness. 
For our main example, we use $n_c=5$ control pairs with driving frequencies $\{u^m_i\}_{i=1,m=1}^{2,n_c}$ in \eqref{eq:control_structure} as they are given in \cref{tab:driving_frequencies_OOP}.
These frequencies are fixed manually within the linewidth (that is $2 d_1$) of the resonance mode in order to gain a sufficiently large signal-to-noise ratio.
Moreover, they should not be too close together to be able to distinguish the response signals by the spectrum analyzer.
Applying the driving frequencies to the real-world experiment, we obtain the laboratory data $\{\bz^{\bu^m}_\star\}_{m=1}^{n_c}$.
To use this data in our numerical optimization, we first subtract the noise $\xi>0$ (defined as the average of the data). 
The actual average noise level $\xi$ of the experiment is between $0.2$ $\mu$V and $0.3$ $\mu$V.
\begin{table}
   \centering
    \begin{tabular}{c|c||c|c}\toprule
       $u^1_1$ & 6.94016   & $u^1_2$ & 6.94036 \\
        \midrule
         $u^2_1$ & 6.94018   & $u^2_2$ & 6.94034\\
         \midrule
         $u^3_1$& 6.94020   & $u^3_2$ & 6.94032\\
         \midrule
%         $f_{1,4}$ &6.940220   & $f_{2,4}$ & 6.940300\\
%        \hline
         $u^4_1$& 6.94024   & $u^4_2$ & 6.94028\\
        \midrule
         $u^5_1$& 6.94025   & $u^5_2$ & 6.94027\\
         \bottomrule
    \end{tabular}
    \caption{Driving frequencies for the $\widetilde{q}_1$ mode in MHz.}
    \label{tab:driving_frequencies_OOP}
\end{table}

We solve the constrained minimization problem \eqref{eq:online_step_minimization} for a fixed $\nu$ using the interior-point method; see, e.g., \cite{NocedalWright2006NumOpt,LuksanMatonoha2004Iteriorpointnonlinearnonconvex}.
\SecondRevisionComment{This is realized using the Matlab built-in function \texttt{fmincon} with appropriate settings.}
Choosing an initial guess $\bp^0 \equiv 0$ we rely for $\bp^\rmref$ on our experience as described in \cref{sec:ParamIdentProblem}. 
In \cref{tab:reconstruction_damping_coupling_2024-07-27_OOP} we report the transformation parameter $\theta$ and the differences $\bar d_1$, $\bar d_2$ calculated from the reference parameter vector $\bp^\rmref$; i.e. $\bar d_1 = d_1 - \bp_1^\rmref$, $\bar d_2 = d_2 - \bp_2^\rmref$.
As tolerance $\mathrm{tol}$, we use machine precision, the algorithm considered nine iterations, and the final regularization $\nu_9$ was $10^6$.

With the optimal parameter $\theta$, we compute the matrices $\sC^m$ (cf. \eqref{eq:matrix_drifted_eigenfrequencies}), calculate the averaged quantities $\langle\lambda\rangle \coloneqq \nicefrac{1}{n_c} \sum_{m=1}^{n_c} \lambda^m$, and analogously $\langle f_{1} \rangle$, $\langle f_{2} \rangle$, and present the average together with the uncertainty interval in \cref{tab:reconstruction_damping_coupling_2024-07-27_OOP}.

The most interesting one for this work is the coupling parameter $\lambda$.
See also the discussion in \cref{sec:Comparison}.

\begin{table}[H]
    \centering
    \begin{tabular}{c|c||c|c||c|c}
    \toprule
        $f_{1}$ & 6.9522 $\pm$ 0.0001 MHz  & 
        $f_{2}$ & 7.0156 $\pm$ 0.0001 MHz & 
        $\lambda$ & 0.6474 $\pm$ 0.0002 MHz\\
        \midrule
        $\theta$ & %1.949802357956
        1.9498 [-]& 
        $\bar d_1$ & 1.5828 $\mu$Hz & 
        $\bar d_2$ & 6.4871 $\mu$Hz\\
        \bottomrule
    \end{tabular}
    \caption{Optimal estimated coefficients. }
    %https://gitlab.inf.uni-konstanz.de/jan.bartsch/oscillators/-/tree/IterativeSchemeConverged
    \label{tab:reconstruction_damping_coupling_2024-07-27_OOP}.
\end{table}

\begin{figure}
	\centering
	\begin{subfigure}[l]{1\textwidth}
		\includegraphics[width=\textwidth]{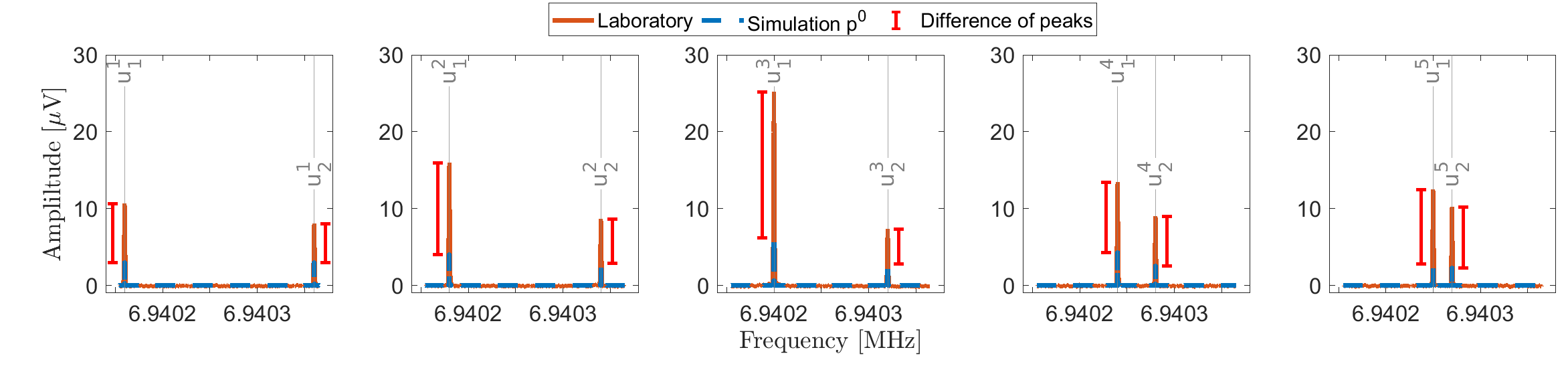}  
        \caption{}
        \label{fig:final_initial_data}
	\end{subfigure}
\\
	\begin{subfigure}[l]{1\textwidth}
		\includegraphics[width=\textwidth]{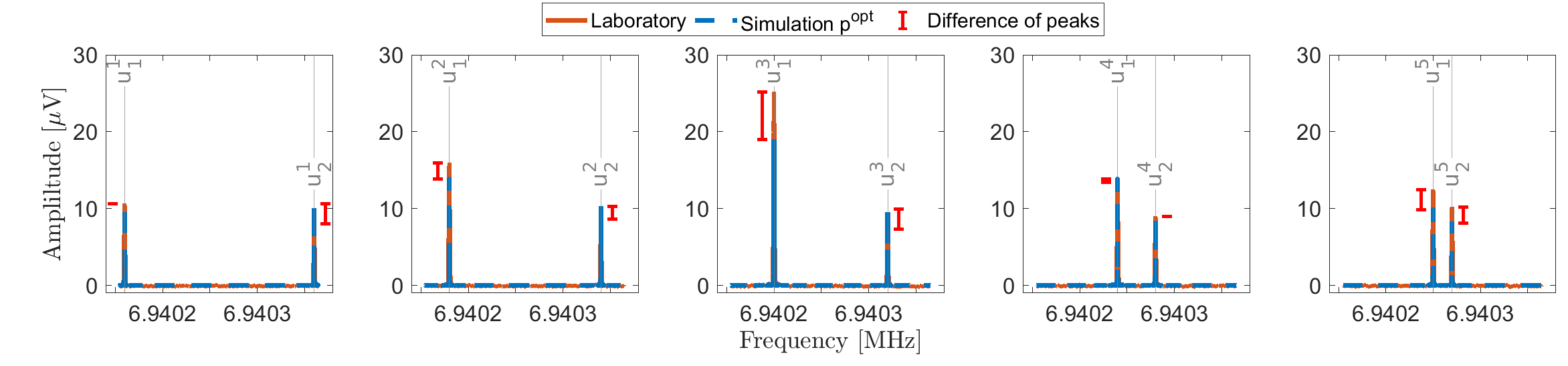}
        \caption{}
        \label{fig:final_recon_data}
	\end{subfigure}
%https://gitlab.inf.uni-konstanz.de/jan.bartsch/oscillators/-/tree/IterativeSchemeConverged
	\caption{Experimental data (solid, orange) and simulated data (dashed, blue) applying the driving frequencies given in \cref{tab:driving_frequencies_OOP}.
    The red errorbars indicate the difference between the amplitudes of the peaks of experimental and simulation data.
     The vertical lines depict the control frequencies;
    (a) before the reconstruction (using $\bp^0\equiv 0$);
    (b) after the reconstruction using $\bp^\rmopt$; 
  for the reconstructed parameters see \cref{tab:reconstruction_damping_coupling_2024-07-27_OOP} and for simulation parameters see \cref{tab:simulation_hyper_parameters_damping_coupling}.}
	\label{fig:results_final_ident_lab}
\end{figure}

The comparison of the laboratory measurements and the simulation results using the initial parameters $\bp^0$ and the optimal ones $\bp^\rmopt$ resulting from \cref{algo:Recon} are presented in \cref{fig:results_final_ident_lab}.
The vertical lines show the driving frequencies.
In \cref{fig:final_initial_data}, we plot the experimental data (solid, orange) and the simulated data (dashed, blue) for the unreconstructed parameters $\bp^0$.
As expected, we observe the two peaks that correspond to the excitation frequencies in all plots.
The peaks have in theory a delta shape since they are Fourier transformations of an essentially sinusoidal signal. 
However, since we work with a finite resolution in simulation and experiment, this approximation results in a small non-zero width of the peak.

With our optimization procedure, a better agreement of the amplitude of the peaks for simulation and experiment is obtained.
In \cref{fig:final_recon_data}, we plot the experimental data and the simulated data with the reconstructed $\bp^\rmopt$ provided by \cref{algo:Recon}.
In fact, in \cref{fig:final_recon_data} we see that the agreement of the simulated and the experimental data is significantly enhanced when using the optimal parameters $\bp^\rmopt$.
This can be observed in particular by the length of the errorbars. 
They demonstrate the difference between the maxima of the peaks in the experimental and simulation data.

\begin{figure}
	\begin{subfigure}[l]{0.45\textwidth}
		\includegraphics[width=\textwidth]{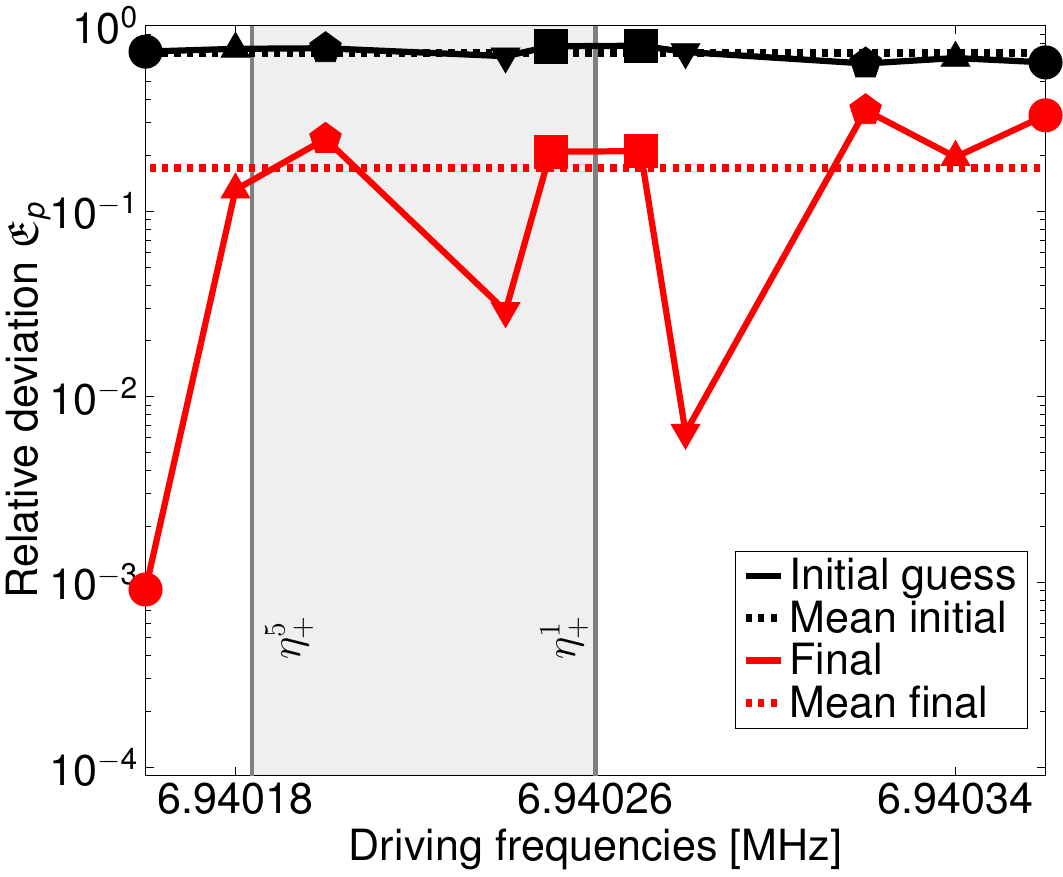}  
		\caption{}
		\label{fig:final_2024-06-26_difference}
	\end{subfigure}
\hfill	
	\begin{subfigure}[l]{0.45\textwidth}
		\includegraphics[width=\textwidth]{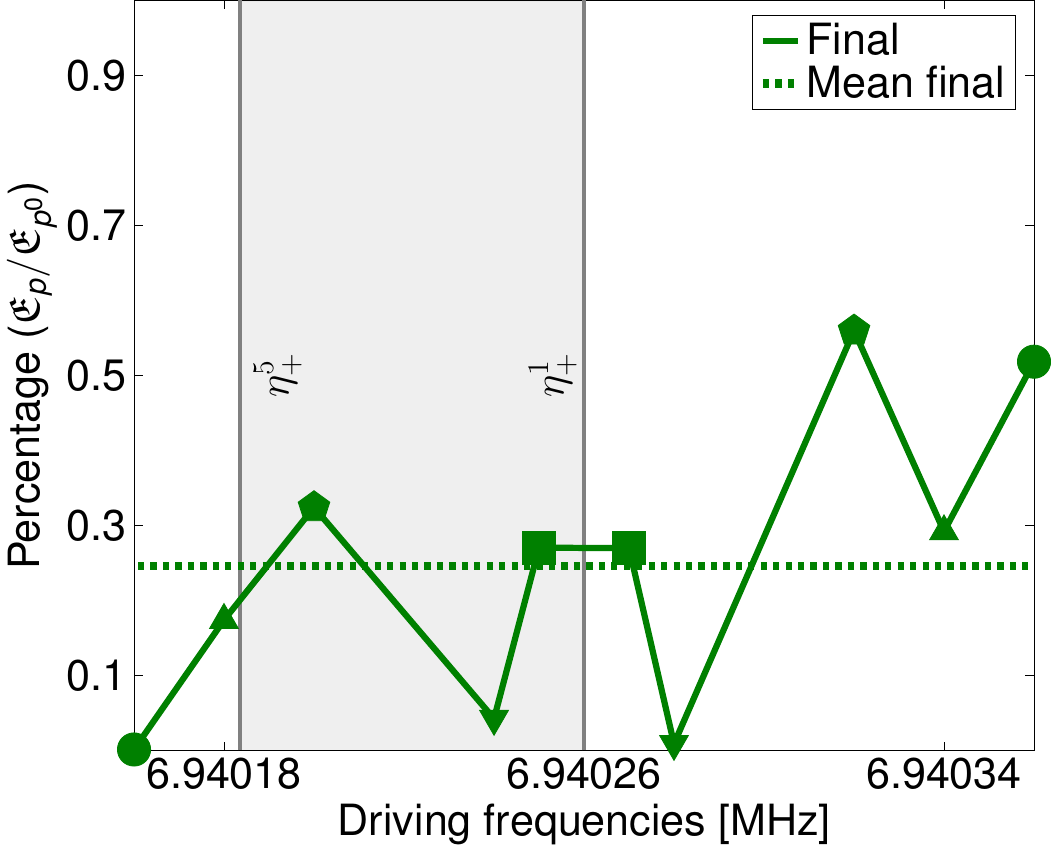}  
		\caption{}
		\label{fig:final_2024-06-26_percentage}
	\end{subfigure}
	\caption{Deviation plots of the result of \cref{algo:Recon}. 
    The different marker shapes depict the different control pairs (cf. \cref{tab:driving_frequencies_OOP});
 (a) Relative deviation for driving frequencies (cf. \cref{def:deviation}) before ($\mathfrak{E}_{\bp^0}$, black) and after ($\mathfrak{E}_{\bp^\rmopt}$, red) the reconstruction process;
 (b) Deviation improvement $\mathfrak{E}_{\bp^\rmopt}/\mathfrak{E}_{\bp^0}$ in percent  (cf. \cref{def:deviation}).}
	\label{fig:final_2024-06-26_Deviation}
\end{figure}

To quantify the improvement, we plot in \cref{fig:final_2024-06-26_difference} the deviation introduced in \cref{def:deviation} for the initial guess (black) and the reconstructed parameter (red) over the driving frequencies.
The gray area visualizes the range of the drift of the eigenfrequencies during the execution of the laboratory experiment (cf. \cref{sec:DriftingEigenfrequencies}).
The initial eigenfrequency is denoted by $\eta_+^1$ and the final one $\eta_+^5$.
The deviations corresponding to the specific control pair are tagged with a distinct marker shape (from $m=1,\ldots,5$: bullet, triangle, pentagon, diamond, square).
These shapes stay the same throughout the paper.
We observe that we can decrease the relative deviation for each control pair.
More in detail, the average relative deviation decreases from 0.7 to 0.17 (dashed lines).
In \cref{fig:final_2024-06-26_percentage}, we plot the relative improvement of the deviation between experimental and simulated data, this is $\mathfrak{E}_{\bp^\rmopt}/\mathfrak{E}_{\bp^0}$ (cf. \cref{def:deviation}), over the driving frequencies.
We see that the deviation between the experimental data and the simulations using $\bp^\rmopt$ is, on average, only 20\% of the deviation between the experimental data and the simulation using $\bp^0$.

\begin{figure}
    \centering
    \begin{subfigure}[l]{0.32\textwidth}
        \centering
    	\includegraphics[scale=\scalingThreeFigures]{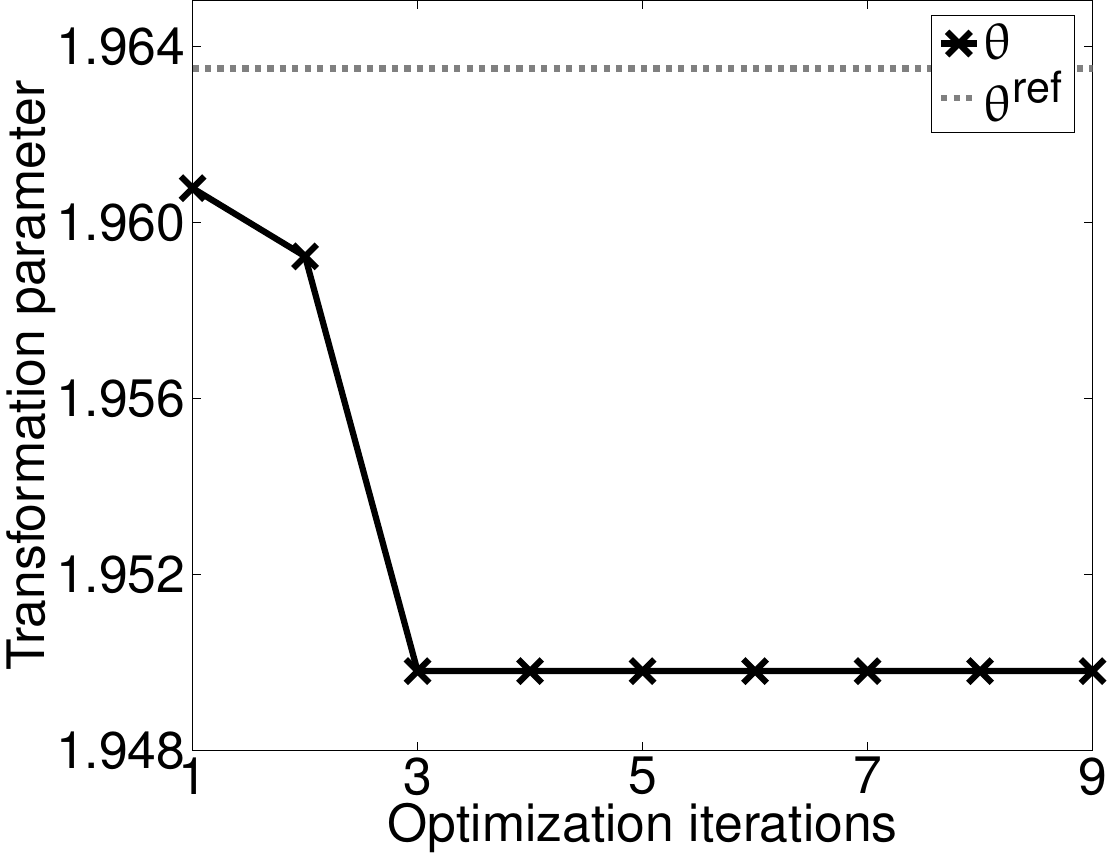}
    	\caption{}
    	\label{fig:conv_theta}
    \end{subfigure}
\hfill
    \begin{subfigure}[l]{0.32\textwidth}
        \centering
        \includegraphics[scale=\scalingThreeFigures]{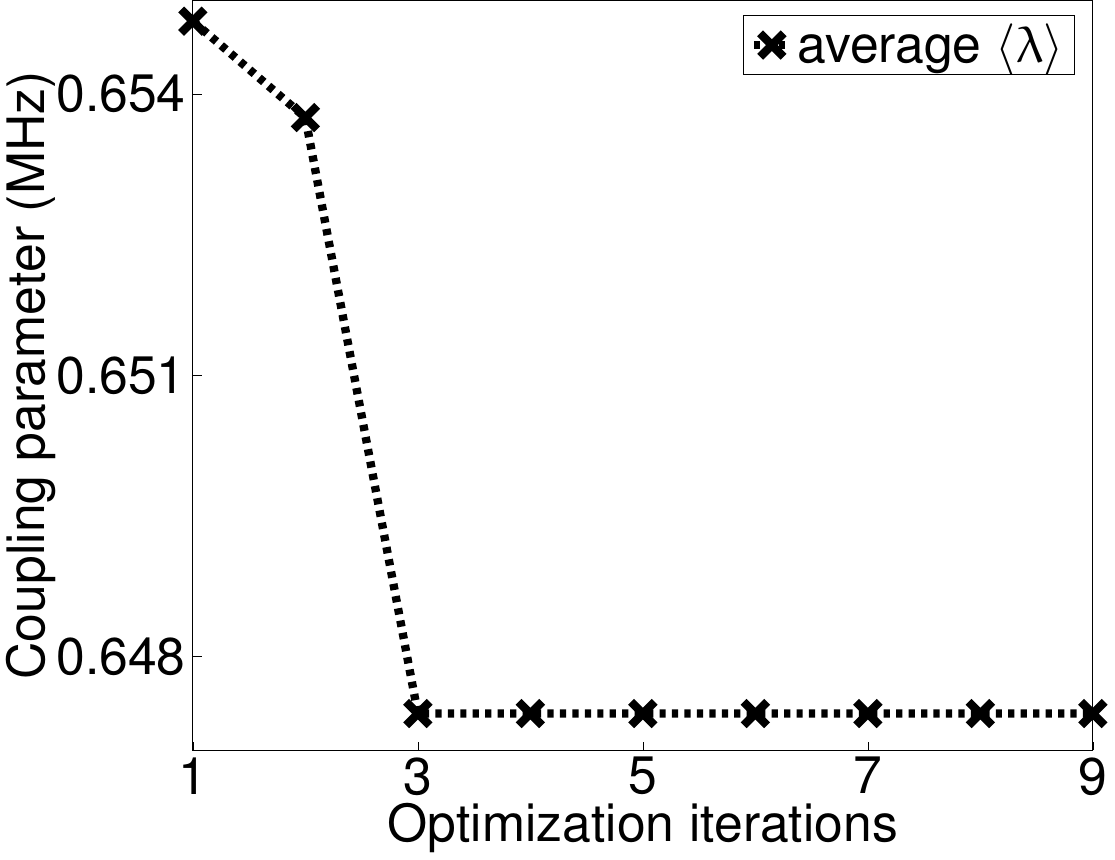}
        \caption{}
        \label{fig:conv_lambda}
    \end{subfigure}
\hfill
    \begin{subfigure}[l]{0.32\textwidth}
        \centering
       \includegraphics[scale=\scalingThreeFigures]{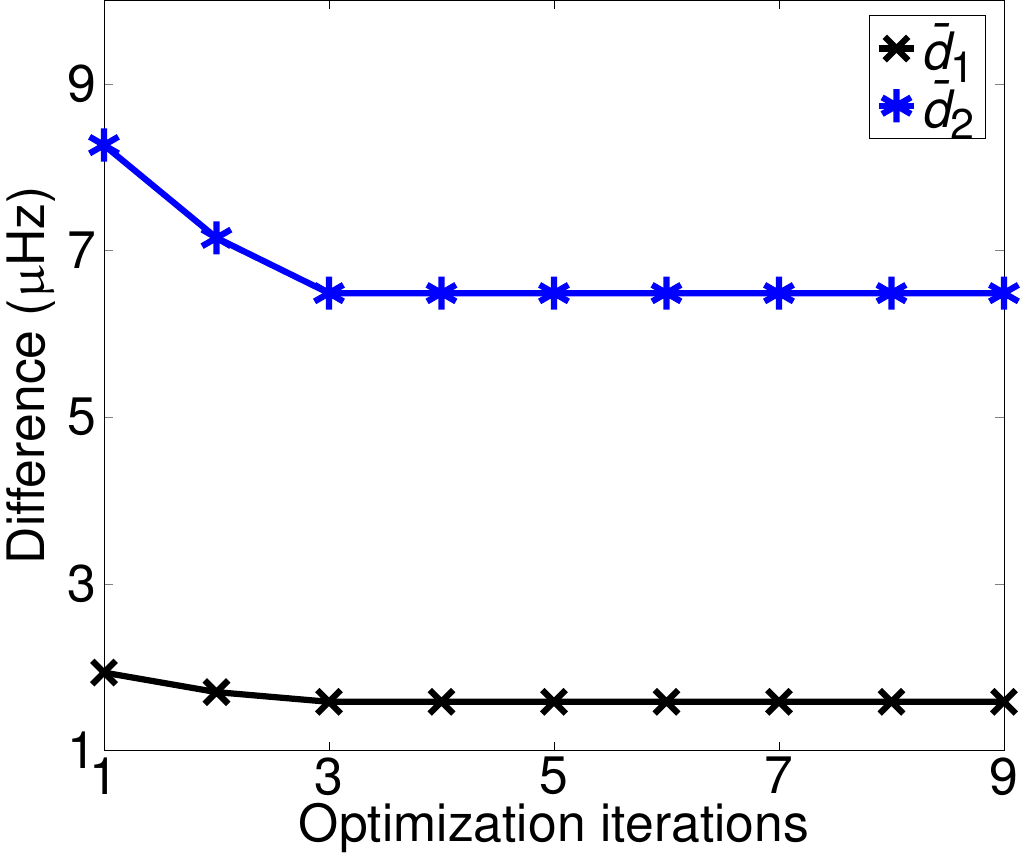}   
        \caption{}
        \label{fig:conv_damping}
    \end{subfigure}
    \caption{Convergence history of \cref{algo:Recon}.
    (a) Convergence of $\theta$;
    (b) Convergence of averaged $\langle\lambda\rangle$;
    (c) Convergence of corrections $\bar{d}_1$, $\bar{d}_2$  to the damping parameters.}
    \label{fig:convergence}
\end{figure}

In \cref{fig:convergence}, we visualize the convergence behavior of \cref{algo:Recon}.
As expected, we see in \cref{fig:conv_theta} that at the first iterations, the value of $\theta$ is close to $\theta^\rmref$ and then converges to a value different from $\theta^\rmref$.
In all three panels, we observe that after the third iteration, the algorithm seems to have converged and only small improvements are made afterwards.

\begin{figure}
	\begin{subfigure}[l]{0.32\textwidth}
		\centering
\includegraphics[scale=\scalingThreeFigures]{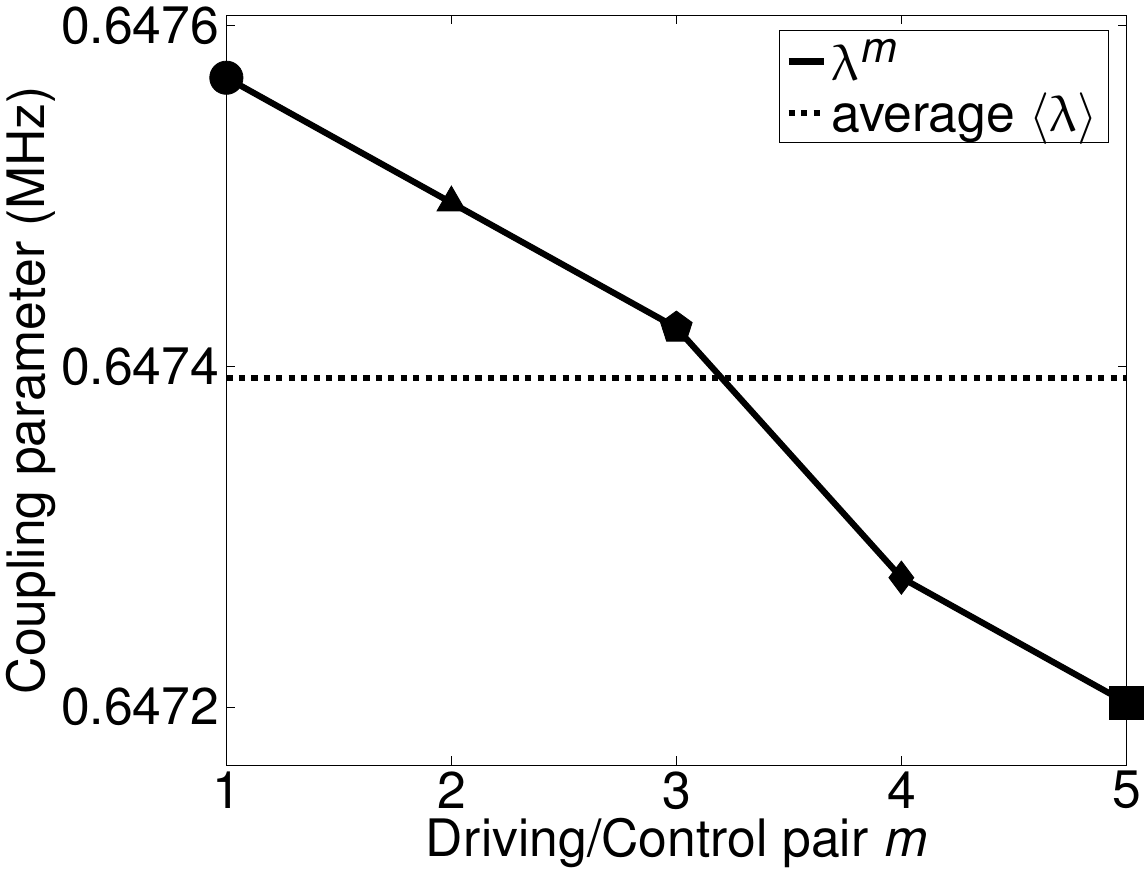}
		\caption{}
		\label{fig:avg_lambda_OOP}
	\end{subfigure}
\hfill
	\begin{subfigure}[l]{0.32\textwidth}
		\centering
		\includegraphics[scale=\scalingThreeFigures]{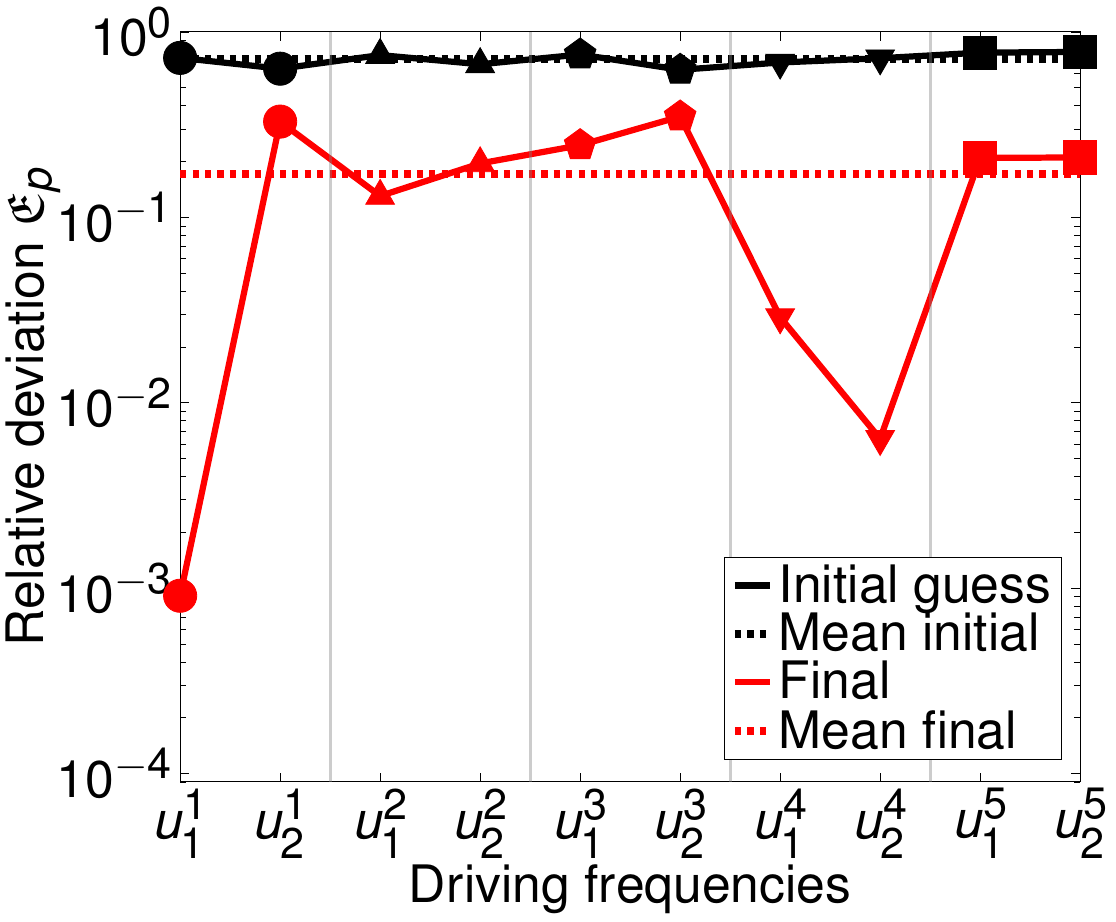}
		\caption{}
		\label{fig:InvidividualSummands_error}
	\end{subfigure}
    \caption{Results of final identification for $\widetilde{q}_1$ mode. The different marker shapes depict the different control pairs given in \cref{tab:driving_frequencies_OOP};
    (a) Coupling constant $\lambda^m$ for controls pairs and average $\langle\lambda\rangle$;
    (b) Deviation for the initial guess ($\mathfrak{E}_{\bp^0}$, black) and in the final iteration ($\mathfrak{E}_{\bp^\rmopt}$, red) ordered by control pairs (cf. \cref{fig:final_2024-06-26_difference}).}
\end{figure}

In \cref{fig:avg_lambda_OOP}, we plot the $\lambda^m, m =1,\ldots,n_c$ and the average $\langle\lambda\rangle$. 
We observe that the difference between the $\lambda^m$ and $\langle\lambda\rangle$ is quite small and in the order of magnitude of the drifting eigenfrequencies.
In \cref{fig:InvidividualSummands_error}, we plot the deviation introduced in \cref{def:deviation} over the control frequency pairs.
We observe that the deviation is higher than the average for the third control pair and the least for the fourth control pair.

\bigskip

To validate the optimal parameters found by \cref{algo:Recon}, we perform a second test.
For this, we choose control frequencies located around the eigenfrequency $\eta_-$ of the mode $\widetilde q_2$.
Then, we solve \eqref{eq:online_step_minimization} once with fixed (small) regularization parameter $\nu$ and use the controls given in \cref{tab:driving_frequencies_IP}.
\begin{table}[H]
	\centering
	\begin{tabular}{c|c||c|c}
    \toprule
		$u^1_1$ &7.02750   & $u^1_2$ & 7.02770\\
		 \midrule
		$u^2_1$& 7.02752   & $u^2_2$ & 7.02768\\
		 \midrule
		$u^3_1$& 7.02754   & $u^3_2$ & 7.02766\\
		 \midrule
		$u^4_1$& 7.02758   & $u^4_2$ & 7.02762\\
		 \midrule
		$u^5_1$& 7.02759   & $u^5_2$ & 7.02761\\
    \bottomrule
	\end{tabular}
	\caption{Driving frequencies around the eigenfrequency $\eta_-$ in MHz.}
	\label{tab:driving_frequencies_IP}
\end{table}

\begin{figure}
    \centering
     \begin{subfigure}[l]{\textwidth}
        \includegraphics[width=\textwidth]{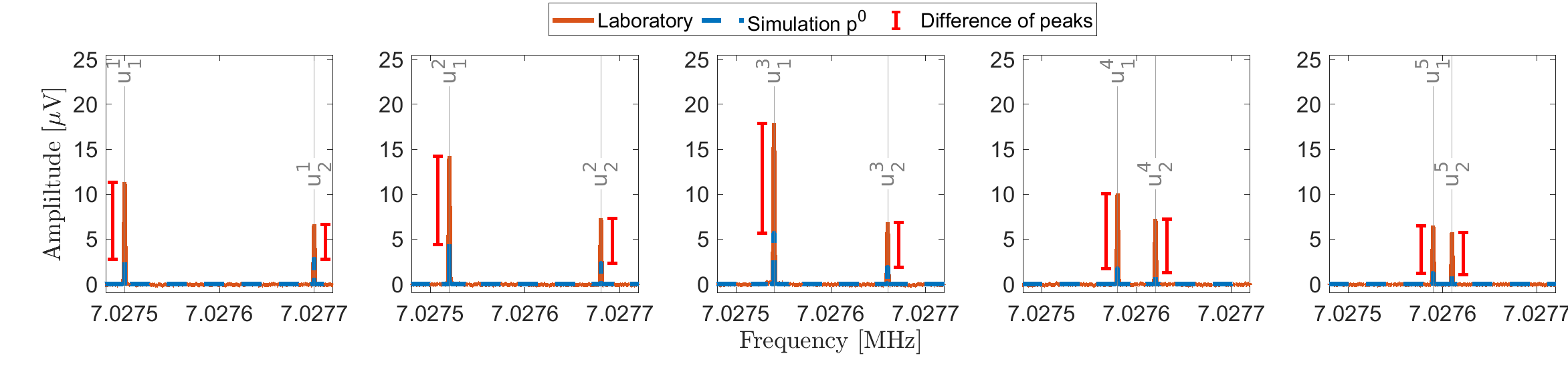}
        \caption{}
        \label{fig:fourier_IP_initial}
    \end{subfigure}
\\
    \begin{subfigure}[l]{\textwidth}
        \includegraphics[width=\textwidth]{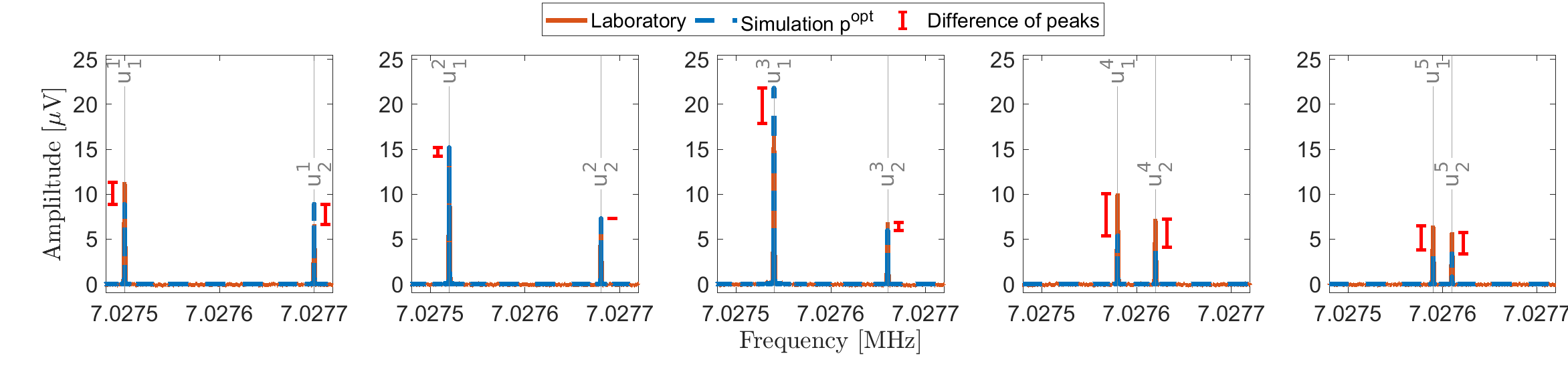}
        \caption{}
        \label{fig:fourier_IP_final}
    \end{subfigure}
    \caption{Experimental data (solid, orange) and simulated data (dashed, blue) applying the driving frequencies given in \cref{tab:driving_frequencies_IP}.
    The red errorbars indicate the difference between the amplitudes of the peaks of experimental and simulation data.
     The vertical lines depict the control frequencies;
    (a) before the reconstruction (using $\bp^0\equiv 0$);
    (b) after the reconstruction using $\bp^\rmopt$; 
     for the reconstructed coefficients see \cref{tab:reconstruction_damping_coupling_2024-07-27_IP} and for simulation parameters \cref{tab:simulation_hyper_parameters_damping_coupling}.}
    \label{fig:fourier_IP}
\end{figure}

In \cref{fig:fourier_IP,fig:error_IP}, we show the results of this optimization.
In \cref{fig:fourier_IP}, we plot similar to \cref{fig:results_final_ident_lab} the experimental data (solid, orange) and the simulated data (dashed, blue) for the unreconstructed parameters $\bp^0$.
The vertical lines show the driving frequencies.
In \cref{fig:fourier_IP_initial} present the Fourier transform for the $\widetilde{q}_2$ mode given the initial guess $\bp^0$ and in \cref{fig:fourier_IP_final} for the optimized parameters $\bp^\rmopt$.
In all plots, the lengths of the error indicate the difference in the amplitude of the peaks in the experimental and simulation data.
\\
In \cref{fig:deviation_IP}, we plot the relative deviation defined in \cref{def:deviation} for the initial guess (black) and the optimized parameters (red).
Furthermore, the gray area depicts the region of the drift of the eigenfrequencies.
Notice that in this experiment, the drift is larger than for the $\tildeq_1$ mode.
In \cref{fig:improvement_IP}, we plot the improvement of the deviation similar to \cref{fig:final_2024-06-26_percentage}. 
Also here we see that we can improve the parameters of the model but not as well as for the $\tildeq_1$ mode.
More specifically, the relative deviation decreases from 0.74 to 0.27 (dashed lines).

Further, we observe that the reconstructed parameters also slightly differ from the ones of the $\widetilde{q}_1$ mode.
Notice that the driving frequencies are farther away from the (shifted) eigenfrequency $\eta^{5}_{-}$ as in the previous experiment.
Hence, we cannot expect a similar or better behavior. 

\def\scalingTwoFigures{0.4}
\begin{figure}
    \centering
     \begin{subfigure}[l]{0.45\textwidth}
        \centering
        \includegraphics[scale=\scalingTwoFigures]{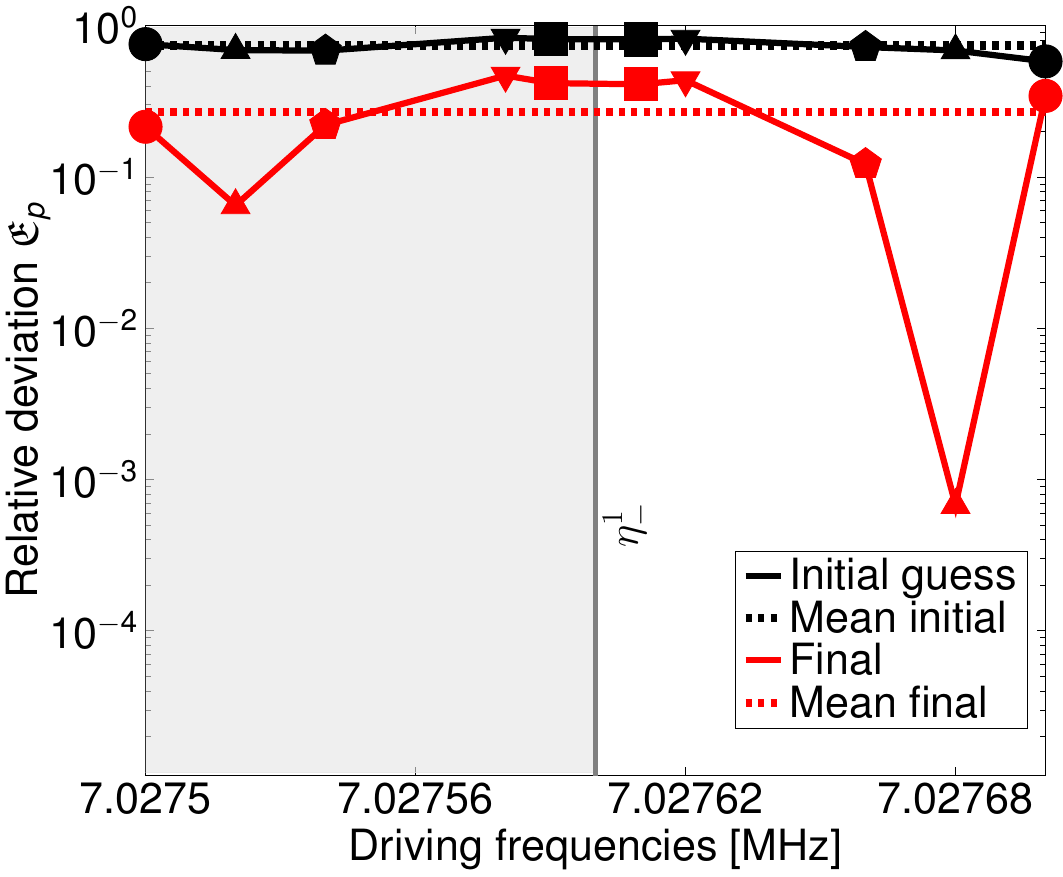}
        \caption{}
        \label{fig:deviation_IP}
    \end{subfigure}
    \hfill
    \begin{subfigure}[l]{0.45\textwidth}
        \centering
        \includegraphics[scale=\scalingTwoFigures]{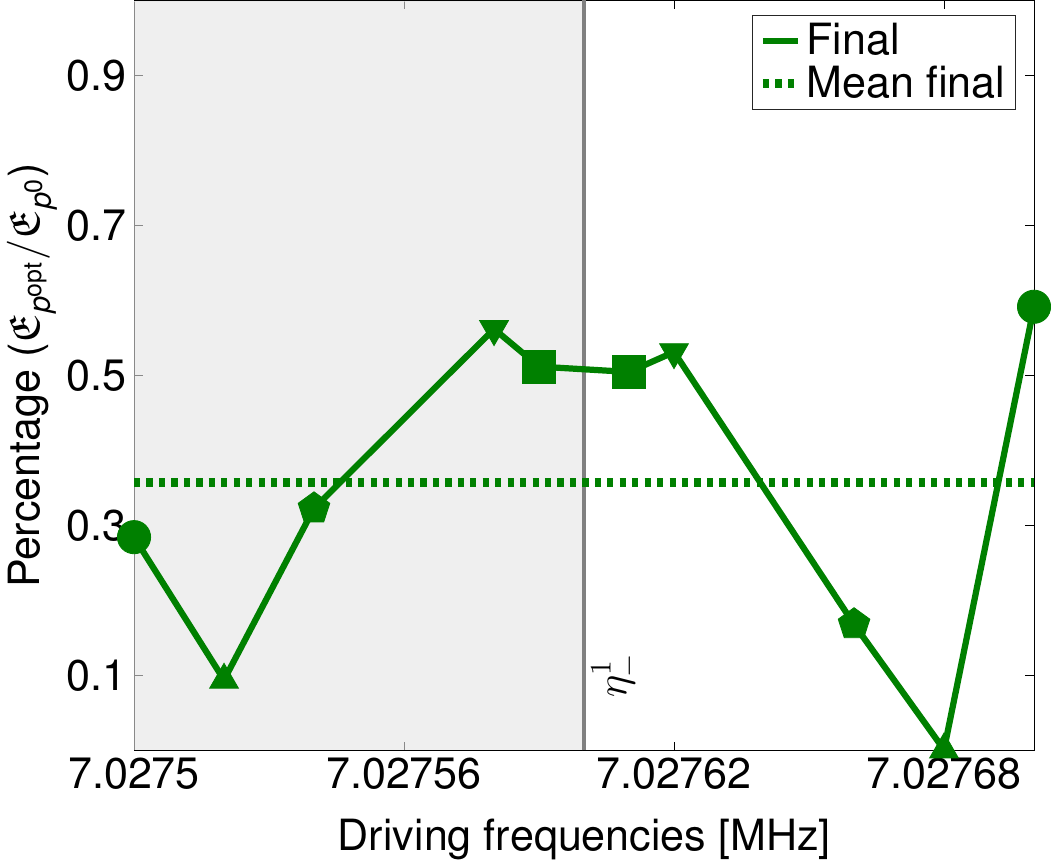}
        \caption{}
        \label{fig:improvement_IP}
    \end{subfigure}
    \caption{Deviation plots and coupling for $\widetilde{q}_2$ mode. The different marker shapes depict the different control pairs (cf. \cref{tab:driving_frequencies_IP});
    (a) Relative deviation $\mathfrak{E}_{\bp}$ (cf. \cref{def:deviation});
    (b) Deviation improvement in percent $\mathfrak{E}_{\bp^\rmopt}/\mathfrak{E}_{\bp^0}$;
    }
    \label{fig:error_IP}
\end{figure}

In \cref{tab:reconstruction_damping_coupling_2024-07-27_IP}, we present the resulting estimated parameters analogous to \cref{tab:reconstruction_damping_coupling_2024-07-27_OOP}.

\begin{table}[H]
	\centering
	\begin{tabular}{c|c||c|c||c|c}
    \toprule
		$f_{1}$ & 6.9582 $\pm$ 0.0001 MHz & 
		$f_{2}$ & 7.009 $\pm$ 0.0001 MHz & 
		$\lambda$ & 0.7016 $\pm$ 0.0002  MHz \\
		\midrule
		$\theta$ & %2.04063418
        2.0406 [-] & 
		$\bar d_1$ &
         6.247 mHz
        %3.92557275e-02 $\frac{\text{rad}}{\text{s}}$
        & 
		$\bar d_2$ & 
        2.906 mHz \\
        %1.82615318e-02  $\frac{\text{rad}}{\text{s}}$\\
    \bottomrule
	\end{tabular}
	\caption{Reconstructed coefficients for $\tildeq_2$ mode. }
	%https://gitlab.inf.uni-konstanz.de/jan.bartsch/oscillators/-/tree/IterativeSchemeConverged
	\label{tab:reconstruction_damping_coupling_2024-07-27_IP}.
\end{table}

%%%%%%%%%%%%%%%%%%%%%%%%%%%%%%%%%%%%%%%%%%%%
\section{Comparison with an independent experiment}
\label{sec:Comparison}
%%%%%%%%%%%%%%%%%%%%%%%%%%%%%%%%%%%%%%%%%%%%
In this section, we compare our results of \cref{sec:NumExp} with those of another experimental approach in \cite{BarakatWeig2024CouplingStrength}. 
In this cited work, the authors used the same string (see Supp. Mat. of \cite{BarakatWeig2024CouplingStrength}) as well as the same experimental setup, thus ensuring the measurement of the same system dynamics.
While in the current work, we use a lock-in amplifier to generate the two excitation signals, the authors in \cite{BarakatWeig2024CouplingStrength} used a noise drive that allows to excite all frequencies within a certain band at the same time. 
Furthermore, an additional harmonic signal was required or a parametric excitation in that experiment.
All other setup details are exactly the same including, in the first place, the investigated sample and the coupled actuation and read-out systems. 
Furthermore, the work in~\cite{BarakatWeig2024CouplingStrength} ranges over an interval of DC voltages $V_{\mathrm{DC}}$ from $-32$ to $32$ V, where the coupling coefficient $\lambda$ is deduced by a \SecondRevisionComment{series of} measurements and analysis based on a mathematical model.

This parallel between the two contributions elevates the certainty in both results since the two approaches are intrinsically different. 
While the cited approach relies on the physical modeling of the system as well as on the use of a coupling-dependent physical phenomenon, called the parametric normal mode splitting (PNMS), to estimate the coupling coefficient \cite{BarakatWeig2024CouplingStrength}, we rely here on a purely mathematical approach, that is, no physical modeling is added in the presented work except for the consideration of a linearly coupled two-mode oscillating system.

Regarding the results, the deduced coupling coefficient from the experiment in~\cite{BarakatWeig2024CouplingStrength} yielded 0.8604 MHz while that of the work presented here gives 0.6474 MHz and 0.7016 MHz for 
%before was 0.7293 MHz found by Ahmed
$\widetilde{q}_1$ and $\widetilde{q}_2$, respectively, nevertheless representing the same coupling coefficient.
The values deduced from the current work and that of \cite{BarakatWeig2024CouplingStrength} are found to be of the same order of magnitude with a difference of about 20\%, which is close to the tolerance given for the numerical experiments inside the current work (see \cref{fig:avg_lambda_OOP}). 
Arriving at these relatively similar results despite using two totally different methodologies, and remarkably despite the usual drifts in experimental results across non-simultaneously performed experiments, supports our findings in this work and thus provides an important asset in extending our methodology to the investigation of other system parameters.

In addition, the work in \cite{BarakatWeig2024CouplingStrength} required a thorough mathematical modeling of the underlying physics, and an additional experiment covering a range of bias voltages to determine the implicit parameters of the system. 
However, in the presented work, we solely require a fast executable experiment at a given bias voltage, without any need for modeling the physical nature of the coupling itself.

In summary, we can state that our approach greatly reduces the experimental overhead in terms of complexity of the experiment, measurement time, number of measurements required, and therefore also the cost of the experiments.

%%%%%%%%%%%%%%%%%%%%%%%%%%%%%%%%%%%%%%%%%%
\section{Conclusion/Discussion}
\label{sec:Conclusion}
%%%%%%%%%%%%%%%%%%%%%%%%%%%%%%%%%%%%%%%%%%

This work presents an efficient approach for reconstructing system coefficients in coupled harmonic oscillators through an iterative optimization method based on Tikhonov regularization and considering actual data from laboratory experiments. 
The success of our method lies in its ability to automatically combine simulation and experimental data, leveraging the structure of the inverse problem. 
By utilizing this approach, our method achieved estimates of the coupling and damping coefficients that coincide very well with a completely different physical approach. 

Our reconstructed parameters are within the correct order of magnitude, indicating that our method reliably captures the system's underlying behavior. 
Given the inherent noise and complexity of real-world experiments, we cannot expect to achieve absolute precision. However, the results show a high level of consistency and alignment with theoretical expectations, which is sufficient for many practical applications.

Another significant achievement of this method is the reduction of laboratory time. 
By employing a combined approach of simulations and minimal experiments, we drastically cut down the number of physical tests required. 
This not only saves time but also reduces the cost associated with running numerous experiments in a laboratory setting.

Importantly, one of the strengths of our method is that it does not require deep knowledge of the underlying physics of the system. 
Instead, the method focuses on reconstructing the system coefficients from observed data, making it broadly applicable to other similar systems. 
This versatility, combined with its efficiency, makes it a valuable tool for system identification and optimization tasks in various scientific and engineering contexts.

One possible improvement of our method is to find optimized driving frequencies selected in a greedy manner. 
This approach has been discussed and applied to small problems in \cite{BuchwaldCiaramella2021GreedyReconstruction}.
By choosing these frequencies carefully, it is possible to capture the necessary dynamics with even less experimentation. 
However, given the current setup, the constraints on the driving frequencies are so stringent that further optimization is unlikely to yield significantly different results.

Furthermore, this approach can also be applied to models in which an unknown nonlinearity has to be reconstructed.
In this case, one can for example choose a set of basis functions that can be used to approximate the nonlinearity.
The optimization problem is then to find the coefficients in a linear combination of the basis elements that lead to the best agreement of simulation and experimental data; see, e.g., \cite{BartschBuchwaldCiaramellaVolkwein2024NonlinRecon2D} for further information on this procedure applied to a semilinear elliptic partial differential equation.

%%%%%%%%%%%%%%%%%%%%%%%%%%%%
\section*{Statements and Declarations}
%%%%%%%%%%%%%%%%%%%%%%%%%%%%

\paragraph*{Funding.}\, This work was funded by the Deutsche Forschungsgemeinschaft within SFB 1432, Project-ID 425217212.

\paragraph*{Acknowledgments}. \revisionComment{We thank the anonymous referees for their helpful questions and remarks.}

\paragraph*{Author contributions.}\, 
All authors worked together in formulating the identification problem, and in writing and editing the manuscript.
J.B. worked on the algorithm to solve the identification problem, implemented the simulation codes, and conducted the final identification.
A.B. contributed to the physical modeling and conducted the laboratory experiments.
S.B. provided contributions to the strategies used in the optimization procedure.
G.C. and S.V. supervised the project and provided conceptualization from a mathematical point of view.
E.W. supervised the project and provided conceptualization from a physical point of view.
G.C., S.V., and E.W. procured funding within the SFB1432.

\paragraph*{Competing interests.}\, The authors have no relevant financial or non-financial interests to disclose.

\paragraph*{Ethics approval and consent to participate.}\, Not applicable.

\paragraph*{Consent for publication.} \, Not applicable.

\paragraph*{Availability of data and materials.} \, The data used in this work can be shared upon reasonable request.

\paragraph*{Acknowledgments.}\,  Not applicable.

%%%%%%
%% References
%%%%%%
\bibliographystyle{acm}
\bibliography{literatureReconOsc}{}

\end{document}